\documentclass[11pt]{article}

\RequirePackage[amsthm,amsmath,natbib]{imsart}

 \parskip        \baselineskip
\usepackage[mathscr]{eucal}
\usepackage{graphicx}
\usepackage{latexsym}
\usepackage{amssymb}
\usepackage{psfrag}
\usepackage{epsfig}
\DeclareMathAlphabet{\mathpzc}{OT1}{pzc}{m}{it}

\include{psbox}

\usepackage{amsfonts}

\usepackage{natbib}
\usepackage{graphicx}

\usepackage{amsfonts}

\usepackage{color}
\usepackage{tikz}
\usepackage[textsize=small]{todonotes}

\begin{document}

\newtheorem{theorem}{\bf Theorem}[section]
\newtheorem{example}{\bf Example}[section]
\newtheorem{definition}{\bf Definition}[section]
\newtheorem{corollary}{\bf Corollary}[section]
\newtheorem{remark}{\bf Remark}[section]
\newtheorem{lemma}{\bf Lemma}[section]
\newtheorem{assumption}{Assumption}[section]
\newtheorem{condition}{\bf Condition}[section]
\newtheorem{proposition}{\bf Proposition}[section]
\newtheorem{definitions}{\bf Definition}[section]
\numberwithin{equation}{section}
\newtheorem{conjecture}{\bf Conjecture}[section]

\newcommand{\be}{\begin{equation}}
\newcommand{\ee}{\end{equation}}
\newcommand{\bes}{\begin{equation*}}
\newcommand{\ees}{\end{equation*}}
\newcommand{\ba}{\begin{aligned}}
\newcommand{\ea}{\end{aligned}}
\newcommand{\bi}{\begin{itemize}}
\newcommand{\ei}{\end{itemize}}
\newcommand{\barq}{\bar{Q}}
\newcommand{\bara}{\bar{A}}
\newcommand{\bard}{\bar{D}}
\newcommand{\bari}{\bar I}
\newcommand{\barw}{\bar{W}}
\newcommand{\barn}{\bar{N}}
\newcommand{\hatq}{\hat{Q}}
\newcommand{\hata}{\hat{A}}
\newcommand{\hatd}{\hat{D}}
\newcommand{\hatj}{\hat{J}}
\newcommand{\hati}{\hat{I}}
\newcommand{\hatw}{\hat{W}}
\newcommand{\hatx}{\hat{X}}
\newcommand{\hatn}{\hat{N}}
\newcommand{\bart}{\bar{T}}
\newcommand{\dott}{\dot{T}}
\newcommand{\tilxi}{\tilde{\xi}}
\newcommand{\tilnu}{\tilde{\theta}}
\newcommand{\tilmu}{\tilde{\mu}}
\newcommand{\tillam}{\tilde{\lambda}}
\newcommand{\tilyn}{\tilde{Y}^n}
\newcommand{\tily}{\tilde{Y}}
\newcommand{\tilc}{\tilde c}
\newcommand{\tilx}{\tilde{X}}
\newcommand{\NN}{\mathbb{N}}
\newcommand{\EE}{\mathbb{E}}
\newcommand{\PP}{\mathbb{P}}
\newcommand{\KK}{\mathbb{K}}
\newcommand{\QQ}{\mathbb{Q}}
\newcommand{\ink}{$\blacksquare$}
\newcommand{\clo}{\mathcal{O}}
\newcommand{\clf}{\mathcal{F}}
\newcommand{\clg}{\mathcal{G}}
\newcommand{\tilclf}{\tilde{\mathcal{F}}}
\newcommand{\LL}{\mathbb{L}}
\newcommand{\RR}{\mathbb{R}}
\newcommand{\qq}{\mathbb{Q}}
\newcommand{\bm}{{B\!M}}
\newcommand{\go}{\rightarrow}
\newcommand{\Go}{\Rightarrow}
\newcommand{\Red}{\textcolor{red}}
\newcommand{\diag}{\mbox{diag}}
\newcommand{\noi}{\noindent}
\newcommand{\skp}{\vspace{\baselineskip}}
\newcommand{\h}{c_K\mu_K^*}
\newcommand{\K}{K_0}
\newcommand{\Y}{\tilde{Y}^n}
\newcommand{\pe}{K}
\newcommand{\T}{\mathbf{T}}
\newcommand{\I}{\mathbf{I}}
\newcommand{\Q}{\mathbf{Q}}
\newcommand{\W}{\mathbf{W}}
\newcommand{\X}{\mathbf{X}}
\newcommand{\hatxx}{\hat{\mathbf{X}}}
\newcommand{\clq}{\mathcal{Q}}
\newcommand{\cla}{\mathcal{A}}
\newcommand{\cld}{\mathcal{D}}
\newcommand{\tilq}{\tilde{Q}}
\newcommand{\tili}{\tilde{I}}
\newcommand{\tileta}{\tilde{\eta}}
\newcommand{\tilw}{\tilde{W}}
\newcommand{\p}{\sigma}
\newcommand{\nt}{\lfloor \varsigma(n) t \rfloor}
\newcommand{\parvs}{\varsigma}
\newcommand{\betan}{\varsigma(n)}
\newcommand{\m}{\mathcal{M}}

\newcommand{\beginsec}{
\setcounter{lemma}{0} \setcounter{theorem}{0}
\setcounter{corollary}{0} \setcounter{definition}{0}
\setcounter{example}{0} \setcounter{proposition}{0}
\setcounter{condition}{0} \setcounter{assumption}{0}
\setcounter{remark}{0} }

\numberwithin{equation}{section} \numberwithin{lemma}{section}

\begin{frontmatter}
\title{Dynamic Scheduling for Markov Modulated Single-server Multiclass Queueing Systems in Heavy Traffic}

 \runtitle{BCP for Markov modulated parallel queues}

\begin{aug}
\author{Amarjit Budhiraja\thanks{This research is partially supported by the National Science Foundation
		(DMS-1004418, DMS-1016441), the Army Research Office (W911NF-10-1-0158) and the US-Israel Binational Science Foundation (Grant 2008466).}, Arka Ghosh and Xin Liu\\ \ \\
}
\end{aug}

\today

\skp

\begin{abstract}
\noi  This paper studies a scheduling control problem for a single-server multiclass queueing network in heavy traffic, operating  in a changing environment. The changing environment is modeled as a finite state Markov process that modulates the arrival and service rates in the system.  Various cases are considered: fast changing environment, fixed environment and slow changing environment. In each of the cases,  using weak convergence analysis,
in particular functional limit theorems for renewal processes and ergodic Markov processes, it is shown that an appropriate \lq\lq{}averaged\rq\rq{} version of the classical $c\mu$-policy (the  priority policy that favors   classes  with higher values of the product of holding cost $c$ and service rate $\mu$) is asymptotically optimal for an infinite horizon discounted cost criterion.

\noi {\bf AMS 2000 subject classifications:} Primary 60K25, 90B22, 90B35, 90B36; secondary 60J70, 60F05.

\noi {\bf Keywords:} Markov modulated queueing network, multiscale queueing systems, heavy traffic, diffusion approximations, scheduling control, scaling limits, asymptotic optimality, $c\mu$ rule, Brownian control problem (BCP).
\end{abstract}

\end{frontmatter}

\section{Introduction}\label{intr-section}
\beginsec
Heavy traffic modeling  has a long history.  For a small sample of some recent works on heavy traffic analysis, we refer the reader to   \cite{harri1, will1, whitt-book, bgv},  a comprehensive list of references can be found in \cite{kushnerbook,whitt-book}. 
The heavy traffic formulation provides tractable approximations for complex queueing systems that capture broad qualitative features of the networks.  While most works in the literature deal with fixed (or internal network-state-dependent) rates, with the advent of modern wireless networks there is an explosion of research on   models for networks where external factors such as meteorological variables (temperature, humidity, etc) can affect the transmission rates to and from the servers (see \cite{buch-kush} and references therein). 
The current paper deals with the simplest such setting of a multiclass network, namely a single server multiclass queueing system  in heavy traffic, where the arrival and services fluctuate according to an environment process.
In the classical constant rate setting this model has been studied in \cite{vanMieghem-cmu} in the conventional heavy traffic regime.
The queueing system consists of a single server which can process $K$ different classes of jobs ($K \ge 1$). 
The $K$ classes  represent different types of jobs (voice traffic, data traffic etc.), and they have different arrival   and service rate functions. We assume that these functions are modulated by a finite state Markov process that represents the background environment. The arrival and service rates satisfy a suitable heavy traffic assumption  which,  loosely speaking, says that the network capacity and service requirements are balanced in the long run. We consider  a scheduling control problem, where the controller decides  (dynamically, at each time point $t \ge 0$)  which class of jobs should the server process so as to minimize an infinite horizon discounted cost function, which involves a linear holding cost per job  per unit time for each of the $K$ different classes.  The jump rates of the modulating Markov process are modeled by a scaling parameter $\nu$. The value  of $\nu$ can vary from negative to positive -- the larger the value of $\nu$, the faster the environment changes -- and we consider three distinct regimes for $\nu$ (see \eqref{scaling-cases}). We show that 
 under different values of $\nu$, with suitably scaling,  queue length processes stabilize leading to different diffusion approximations. 
Main result of the paper shows that, in each of these regimes, an asymptotically optimal scheduling policy is a variation of the classical and intuitively appealing $c\mu$-rule. 
Classical $c\mu$-rule \cite{orig-cmu} is a simple priority policy, where the priority is always given to the class (of nonempty queue) with the highest $c\mu$-value where $c$ and $\mu$ represent  the holding cost and the service rate of that class, respectively.
It is well known that $c\mu$-rule is optimal in many cases, see \cite{stolyar-gen-cmu, vanMieghem-cmu} and references therein.
In our model, the arrival rates change randomly according to the background Markov process, and it is far from clear if the classical $c\mu$-rule would be optimal. In addition, in one of the regimes we consider, the service rates 
are allowed to be modulated by the random environment as well, which means that the order of the $c\mu$-values  keeps changing according to the environment. In particular, in situations where the environment process is not directly observable and one only knows (or can estimate)
its statistical properties, the classical $c\mu$-rule cannot be implemented. 
 We  propose a modified $c\mu^*$-rule, where $\mu^*$ represents  the average service rate of each class  and the average is  taken with respect to  the stationary distribution of the environment process. Under the $c\mu^*$-rule, the server  always processes jobs from the nonempty queue with the highest $c\mu^*$-value. We show that, under an appropriate heavy traffic condition and a suitable scaling, the $c\mu^*$-rule is asymptotically optimal for the chosen cost functional for this environment-dependent queueing model, in each of the three regimes for $\nu$.

\begin{figure}[t]
\begin{center}
\hspace{1.5cm}   \includegraphics[trim=0in 5in 0in 0in,clip=true, angle=90, width=3in]{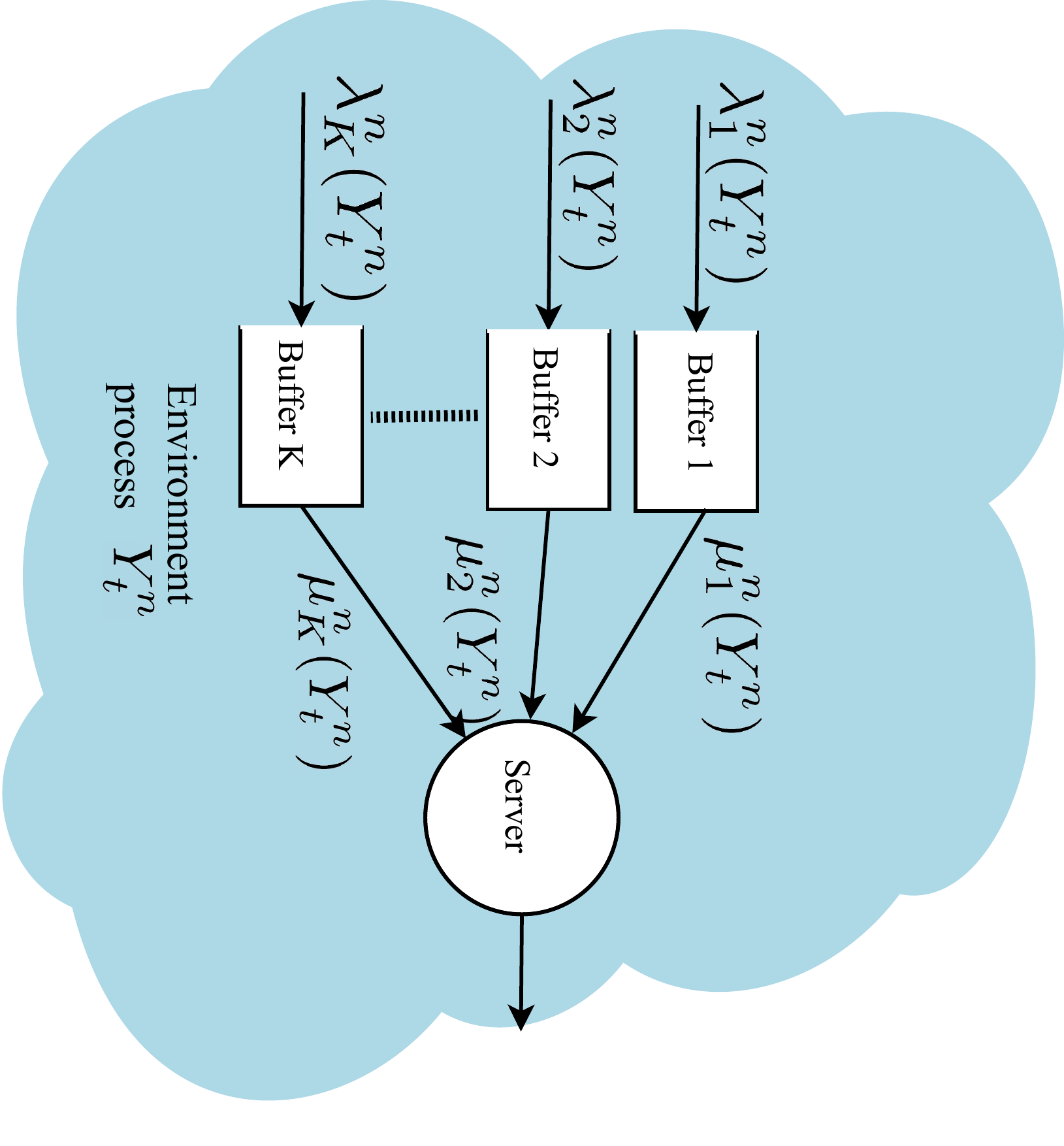}\\
\caption{The $n^{th}$ Markov modulated single-server multiclass queueing system.}\label{fig1}
\end{center}
\end{figure}

As is common practice in heavy traffic analysis, we consider a sequence of networks, indexed by $n \ge 1$, each of which is a single server multiclass network in a random environment (see Fig.~\ref{fig1}). In our formulation, the environment process  operates at a time-scale of $n^\nu$, where  the scaling parameter $\nu$ determines the rate of fluctuations of the environment process.
When $\nu>0$ (Case 1), the background environment is changing faster than  typical arrivals and services (which are $\mathcal{O}(1)$), while the reverse is true when $-1<\nu<0$ (Case 3). If $\nu=0$ (Case 2), the variation rates of the background process and the arrival and service processes are both $\mathcal{O}(1)$. 
The scaling that stabilizes the processes is different in each of the three cases:  When $\nu\geq 0$, the space is scaling down by the usual factor of $n^{\frac{1}{2}}$, and when $1<\nu<0$, the scaling factor becomes $n^{\frac{1-\nu}{2}}$ (time is accelerated by a factor of $n$ in all cases). In each of the cases, we show, using functional central limit theorems for renewal processes and ergodic Markov processes, that the `driving process' converges to a Brownian motion (the variance parameter of the  Brownian motion is different in the three regimes, c.f. Lemma \ref{conv}). 
The application of the first limit theorem is quite common in heavy traffic theory, however the use of limit theorems for ergodic Markov processes in optimal
scheduling problems is less standard. 
We then formulate a Brownian control problem (BCP), which is a formal approximation of the queueing network control problem, where the driving process is replaced by a Brownian motion. 
This BCP is explicitly solvable and  an optimal solution of the BCP is such that,  the optimal state process (formal approximation of the scaled queue-length process) is zero for all the classes except for the one with the lowest value of the corresponding $c\mu^*$ quantity. We  use this insight gained from the solution of the BCP  to propose a priority policy (see Definition \ref{optimal-policy}), which prioritizes queues based on the $c\mu^*$ value of the class. Our main result is Theorem \ref{optimal}, which shows that, under conditions, this   $c\mu^*$  priority policy is  asymptotically optimal for the queueing control problem, in each of the three regimes. 

As noted above, we consider three different regimes for $\nu$: (1) $\nu > 0$; (2) $\nu=0$; (3) $\nu \in (-1,0)$.
Additionally, Case 1 is further subdivided into: (a) $\nu > 1/2$; and (b) $\nu \in (0, 1/2]$.
For cases 1(b), 2 and 3 we take the service rates to be constant (i.e. only the arrival rates are environment dependent).
For case 1(a), both the service and arrival rates are allowed to be environment dependent.
The key result that allows the treatment of environment dependent service rates is Lemma \ref{ergo}, where the condition $\nu > 1/2$ is crucially used.
Although we leave the case when $\nu \le 1/2$ and service rates are environment dependent as an open problem, numerical experiments in Section \ref{main-result}
suggest that even in these settings $c\mu^*$ priority policy may be suitable.

The paper is organized as follows. Section \ref{model} describes the model, introduces the key  assumptions and presents the main scheduling control problem. In Section \ref{main-result}, we introduce the cost criterion and state the main result of the paper, namely Theorem \ref{optimal} (see also Theorem \ref{optimal-limit}). Section \ref{bcp} introduces the Brownian control problem associated with the scheduling control problem from
Section \ref{main-result}. Section \ref{proof}  is devoted to the proofs of Theorems \ref{optimal} and \ref{optimal-limit}.
 Finally, the appendix collects proofs of some auxiliary results and findings of our numerical experiments.

The following notation will be used. Denote the set of positive integers by $\NN$. For $K\in\NN$, let $\RR^K$ be the $K$-dimensional Euclidean space. Define $\mathbb{R}_{>0}^{\otimes K}=\{(x_1, \ldots, x_K): x_i>0, i =1,2,\ldots,K\}$. Vectors are understood to be column vectors. For a $K$-dimensional vector $x$, denote by $\diag(x)$ the $K \times K$ diagonal matrix whose diagonal entries are given by the components of $x$. For $i = 1, \ldots , K$, $x_i$ will denote the $i^{th}$ component of $x$. The transpose of $x$ is denoted by $x'$, and the inner product of two $K$-dimensional vectors $x$ and $y$ will be denoted by $x \cdot y$. Let $|\cdot|$ be a norm on $\RR^K$ given by $|x|=\sqrt{\sum_{i=1}^K x_i^2}$. For a matrix $M$, denote by $M_{ij}$ the $(i,j)^{th}$ entry of $M$. Let
$\mathcal{D}([0,\infty): \RR^K)$ denote the space of functions that are right continuous  with left limits (RCLL) defined from $[0, \infty)$ to $\RR^K$ with the usual Skorohod topology, and $\mathcal{C}([0,\infty), \RR^K)$ the class of continuous functions from $[0, \infty)$ to $\RR^K$ with the local uniform topology. Convergence of random variables $X^n$ to $X$ in distribution will be denoted as $X^n\Rightarrow X$, and weak convergence of probability measures $\mu^n$ to $\mu$ will also be denoted as $\mu^n \Rightarrow \mu$. A sequence of $\mathcal{D}([0,\infty): \RR^K)$-valued random variables $\{X^n\}$ is said to be $C$-tight if and only if the measures induced by $\{X^n\}$ on $(\mathcal{D}([0,\infty): \RR^K), \mathcal{B}(\mathcal{D}([0,\infty): \RR^K)))$ form a tight sequence and any weak limit of the sequence is supported on $\mathcal{C}([0,\infty): \RR^K)$.
For a stochastic process $Y$, we will use the notation $Y(t)$ and $Y_t$ interchangeably, and we let $\Delta Y(t) = Y(t) - Y(t-), \ t\geq 0.$ For semimartingales $X, Y\in \mathcal{D}([0,\infty): \RR)$, we denote by $[X,Y]$ and $\langle X, Y \rangle$ the quadratic covariation and predictable quadratic covariation of $X$ and $Y$, respectively. For a metric space $U$, let $\mathcal{P}(U)$ be the collection of all probability measures on $U$ and $\bm(U)$ the class of bounded measurable functions from $U$ to $\RR$. For $\mu\in \mathcal{P}(U)$ and $f\in \bm(U)$, let 
$
\mu(f) = \int_{U}f(x)\mu(dx).
$
For $f=(f_1,\ldots, f_L)'$, where $f_i\in\bm(U), i=1,2,\ldots,L$, we write
$$
\mu(f) =  \left(\mu(f_1), \ldots, \mu(f_L)\right)'.
$$
A countable set will be endowed with the discrete topology. Finally, we denote generic positive constants by $a_1, a_2, \ldots$, the values of which may change from one proof to another.

\section{Single-server multiclass queueing system}\label{model}
\beginsec
\subsection{Network model}\label{netmod}
Consider a queueing system consisting of a single server which is shared by $K\in\NN$ parallel classes of customers as depicted in Figure~\ref{fig1}.  Consider a sequence of such systems indexed by $n\in\NN$. Let $(\Omega^n, \mathcal{F}^n, \PP^n, \{\mathcal{F}_t^n\}_{t\geq 0})$ be a filtered probability space satisfying the usual conditions. All the random variables and stochastic processes for the $n^{th}$ system are assumed to defined on this space. For the $n^{th}$ system, we define a $\{\mathcal{F}_t^n\}$-Markov process $\{Y^n_t := \Y(n^\nu t): t\geq 0\}$, where $\nu\in\RR$ and $\tilde{Y}^n$ is a $\{\mathcal{F}_{t/n^\nu}^n\}$ Markov process with state space $\LL:=\{1, 2, \ldots, L\}$, a unique stationary distribution $\pi^n:=\{\pi_1^n, \pi_2^n,\ldots, \pi_L^n\}$, and infinitesimal generator (rate matrix) $\qq^n$ which converges to the rate matrix $\qq$
of an ergodic Markov process.  Note that as $n \to \infty$, $\pi^n$ converges to the unique stationary distribution $\pi$ of this Markov process.  
The expectation operator with respect to $\PP^n$ will be denoted by $\EE^n$, but frequently we will suppress $n$ from the notation.
We will make the following uniform ergodicity assumption.
\begin{assumption}\label{uniferg}
	For some $R \in (0, \infty)$ and $\rho \in (0,1)$
$$\sup_n \max_{i,j \in \LL}| P^{n,t}(i,j)- \pi^n_j| \le R \rho^t, \; t \ge 0,$$	
where $ P^{n,t}(i,j) = \PP^n\left(\tilde Y^n(t) = j \mid \tilde Y^n(0) = i\right)$.
\end{assumption}
As shown in Figure~\ref{fig1}, there are $K$ classes of external arrivals. For $i\in\KK:=\{1, 2,\ldots,K\}$, customers of Class $i$ arrive from outside to buffer $i$ and upon completion of service, they leave the system. For $t\geq 0$, denote by $A^n_i(t)$ the number of customers of Class $i$ arriving from the outside by time $t$ and $S^n_i(t)$ the number of customers of Class $i$ completing service by time $t$ if the server is continuously working on Class $i$ customers during $[0,t]$. We assume $A^n_i$ and $S^n_i$ are Poisson processes with intensities modulated by Markov process $Y^n$. More precisely, for $i\in\KK$, let
$\lambda^n_i, \mu^n_i: \LL\go\RR_+$. Then $A^n_i$ and $ S^n_i, i\in\KK$ can be described as follows: For $t\geq 0$,
\be \label{eq:eq5.1}
A^n_i(t)= N_{a,i}^n\left(\int_0^t \lambda^n_i(Y^n_s)ds\right), \;  S^n_i(t)= N_{s,i}^n\left(\int_0^t \mu^n_i(Y^n_s)ds\right),
\ee
where $N_{a,i}^n$ and $N_{s,i}^n$ are independent unit rate Poisson processes. Here $\lambda^n_i, \mu^n_j$ are the  arrival and service rate functions, respectively. 
The processes $Y^n$ and $N^n_{a,i}, N^n_{s,j}, i, j\in\KK$ are assumed to be mutually independent. Also, we assume that
$A_i^n(t) - \int_0^t \lambda_i^n(Y^n(s)) ds$ is a $\{\clf^n_t\}_{t\ge 0}$ martingale.  An analogous assumption on $S_i^n$ will be introduced in Section \ref{scsp} (see Condition \ref{policy-cond} (iii)).
Let $\lambda^n=(\lambda^n_1,\ldots, \lambda^n_K)'$ and $\mu^n=(\mu^n_1,\ldots, \mu^n_K)'$.
We make the following convergence assumptions on $\lambda^n, \mu^n$. Define  the ``averaged" arrival and service rates as follows: 
\bes
\lambda^{n,*} = \pi^n(\lambda^n), \;\;  \mu^{n,*}=\pi^n(\mu^n).
\ees
\begin{assumption}\label{assump}\hfill
\bi
\item[\rm (i)] There exists a  function $\lambda: \LL\go \mathbb{R}_{>0}^{\otimes K}$ such that as $n\go\infty$,
$
\lambda^{n}(y)\go\lambda(y)$, for all $y\in\LL$. 
\item[\rm (ii)] There exists $\mu^*\in\mathbb{R}_{>0}^{\otimes K}$ such that as $n\go\infty$,
$\mu^{n,*}\go \mu^*$.
\ei
\end{assumption}
Let $\lambda^* = \pi(\lambda)$.
Assumption (i) in particular says that the arrival rates are $\clo(1)$ and  $\lim_{n\to\infty}\lambda^{n,*}=\lambda^*.$
 Assumption on the service rates is somewhat weaker in that it only requires that the  ``averaged" service rates are convergent. Additional conditions on $\lambda^{n}$ and $\mu^{n}$ (heavy traffic condition) are formulated in Assumption \ref{ht-assump} at the end of Section \ref{scsp}.

\subsection{Scheduling control and scaled processes}\label{scsp}
In the $n^{th}$ network, the scheduling policy is described by a service allocation process 
$$T^n(t) = (T^n_1(t), T^n_2(t), \ldots, T^n_K(t))', t\geq 0,$$ where for $i\in\KK$, $T^n_i(t)$ denotes the total time used to serve customers of Class $i$ during $[0,t]$. The idle time process is defined as follows:
\be\label{idle}
I^n(t) = t - \sum_{i=1}^KT^n_i(t), \;t\geq 0,
\ee
which gives the total idle time of the server during $[0,t]$. 
Denote by $Q^n$ the $K$-dimensional state process, i.e., the queue length process (including the job being served). Then the evolution of the queueing system can be described as follows: For $i\in\KK$, 
\bes
Q^n_i(t) = N_{a,i}^n\left(\int_0^t \lambda^n_i(Y^n_s)ds\right) - N_{s,i}^n\left(\int_0^t \mu^n_i(Y^n_s)dT^n_i(s)\right), \; t\geq 0.
\ees
Here, for simplicity we assume that at time $0$ the system is empty. 

The processes $T^n$ and $I^n$ are required to satisfy the following conditions: For $i\in\KK$,
\be\label{adm-cond}\ba
&T^n_i(\cdot) \; \mbox{is continuous and nondecreasing, and} \; T^n_i(0)=0,\\
&I^n(\cdot) \; \mbox{is continuous and nondecreasing, and} \; I^n(0)=0,\\
& Q^n_i(t)\geq 0,\; \mbox{for all} \; t\geq 0.
\ea\ee
As a consequence of the  conditions in \eqref{adm-cond}, we have the following Lipschitz  property of $T^n$: For $i\in\KK$,
\be\label{lip}
T^n_i(\cdot) \; \mbox{is uniformly Lipschitz continuous with Lipschitz constant bounded by} \; 1.
\ee
Therefore, $T^n_i(\cdot), i\in\KK$, is almost everywhere differentiable. Denoting its derivative by $\dott^n_i(\cdot)$, we have
\be\label{busy-time}
T^n_i(t) = \int_0^t \dott^n_i(s)ds, \; t\geq 0.
\ee
We  make the following additional assumptions on $\{T^n\}$. Define a sequence of nondecreasing random times $\{\tau^n_k\}_{k\in\NN\cup\{0\}}$.  Let $\tau^n_0=0$, and for $k\in\NN$, $\tau^n_{k+1}$ be the first time after $\tau^n_k$ when there is either an arrival or a service completion. 
\begin{condition}\label{policy-cond}
\item[\rm (i)] For all $k\in\NN\cup\{0\}, \dott^n(t)=\dott^n(\tau^n_k)$ for all $t\in[\tau^n_k, \tau^n_{k+1}).$ 
\item[\rm (ii)] The process $\{\dott^n(t)\}_{t\ge 0}$ is $\{\mathcal{F}^n_t\}$-progressively measurable. 
\item[\rm (iii)] $N_{s,i}^n\left(\int_0^t \mu^n_i(Y^n_s)dT^n_i(s)\right) - \int_0^t \mu^n_i(Y^n_s)dT^n_i(s)$ is a $\{\clf^n_t\}_{t\ge 0}$ martingale.
\end{condition}
Part (i) of the condition says $\dott^n$ does not change values between two successive changes 
of the state of the system.  Part (ii) is a natural nonanticipativity property of $T^n$ and the martingale condition in part (iii) is analogous to the one imposed on $A_i^n$
in Section \ref{netmod}.  We call $T^n$ satisfying \eqref{adm-cond} and Condition \ref{policy-cond} an admissible  control policy and denote the collection of all such $T^n$ by $\cla^n$. 

Using \eqref{busy-time}, we can  write $Q^n_i, i\in\KK$ in the following way: For $t\geq 0$,
\be\label{queue}
Q^n_i(t) = N_{a,i}^n\left(\int_0^t \lambda^n_i(Y^n_s)ds\right) - N_{s,i}^n\left(\int_0^t \mu^n_i(Y^n_s)\dott^n_i(s)ds\right).
\ee
Finally, we write $N^n_a=(N^n_{a,1},\ldots, N^n_{a,K})'$ and $ N^n_s=(N^n_{s,1},\ldots, N^n_{s,K})'$.

We now define the fluid and diffusion scaled processes. Roughly speaking,  the fluid-scaled processes are defined by accelerating time by $n$ and scaling down space  by the same factor, and the diffusion-scaled processes are defined by accelerating time by $n$ and scaling down space
(after an appropriate centering)  by $n^\alpha$ for some $\alpha\in(0,1)$. Choice of $\alpha$ is specified in \eqref{scaling-cases}. 

{\bf Fluid scaling:} For $t\geq 0$, 
\be\label{fluid-scale}\ba
&\barq^n(t) = \frac{Q^n(nt)}{n}, \; \barn^n_{s}(t) = \frac{N^n_{s}(nt)}{n}, \;\barn^n_{a}(t) = \frac{N^n_{a}(nt)}{n}, \\
& \bart^n(t)=\frac{T^n(nt)}{n},\; \bari^n(t) = \frac{I^n(nt)}{n}. 
\ea\ee

{\bf Diffusion scaling:} For $t\geq 0$, 
\be\label{diff-scale}\ba
&\hatq^n(t) = \frac{Q^n(nt)}{n^\alpha}, \; \hatn^n_{a}(t) = \frac{N^n_{a}(nt)-nt}{n^\alpha}, \; \hatn^n_{s}(t) = \frac{N^n_{s}(nt)-nt}{n^\alpha}, \; \hati^n(t)=\frac{I^n(nt)}{n^\alpha}. 
\ea\ee
From the above scalings and \eqref{queue} we see, for $i\in\KK$ and $t\geq 0$,
\be\label{queue-fluid}
\barq^n_i(t)= \barn^n_{a,i}\left(\frac{1}{n}\int_0^{nt}\lambda^n_i(Y^n_s)ds\right)-\barn^n_{s,i}\left(\frac{1}{n}\int_0^{nt}\mu^n_i(Y^n_s)\dott^n_i(s)ds\right),
\ee
and
\be\label{queue-hat1}\ba
\hatq^n_i(t)&=\hatn^n_{a,i}\left(\frac{1}{n}\int_0^{nt}\lambda^n_i(Y^n_s)ds\right)-\hatn^n_{s,i}\left(\frac{1}{n}\int_0^{nt}\mu^n_i(Y^n_s)\dott^n_i(s)ds\right)\\
&\quad+ \frac{1}{n^\alpha}\int_0^{nt}\left(\lambda^n_i(Y^n_s)-\mu^n_i(Y^n_s)\dot{T}^n_i(s)\right)ds.
\ea\ee

We now present our main heavy traffic condition.
\begin{assumption}[Heavy traffic condition]\label{ht-assump}
Let $\lambda^*$ and $\mu^*$ be as in Assumption \ref{assump}. Then
\be\label{hr1}
\sum_{i=1}^K\frac{\lambda_i^*}{\mu^*_i}=1,
\ee
and furthermore, for each $i\in\KK$, there exists $b_i\in\RR$ such that as $n\go \infty$,
\be\label{hr2}
n^{1-\alpha}\left(\frac{\lambda_i^{n,*}}{\mu^{n,*}_i}-\frac{\lambda_i^*}{\mu^*_i}\right)\go b_i.
\ee
\end{assumption}

\section{Objective and main results}\label{main-result}
For our optimization criterion we consider an expected infinite horizon discounted linear holding cost. More precisely, for the $n^{th}$ network, the cost function associated with the control policy $T^n$  is defined as follows:
\be \label{eq:eq1050}
\hat{J}^n(T^n)=\EE\left(\int_0^\infty e^{-\gamma t} (c\cdot\hatq^n(t)) dt\right),
\ee
where $\gamma\in(0,\infty)$ is the ``discount factor" and $c:=(c_1, c_2, \ldots, c_K)\in \mathbb{R}_{>0}^{\otimes K}$ is the vector of ``holding costs" for the $K$ buffers.

The goal of this work is to find a sequence of admissible  control policies $\{\tilde{T}^n:n\in\NN\}$ such that it achieves asymptotic optimality in the sense that
\be\label{asymp}
\lim_{n\go\infty} \hat{J}^n(\tilde{T}^n) = \inf\liminf_{n\go\infty}\hat{J}^n(T^n),
\ee
where the infimum is taken over all admissible control policies $\{T^n\}$.

For single server multiclass queueing systems (with constant rates), one attractive policy is the well-known $c\mu$ rule, where $c$ and $\mu$ stand for holding cost and  service rate, respectively. Under $c\mu$ rule, the server always selects jobs from the nonempty queue with the largest $c\mu$ values. In this work, we generalize the $c\mu$ rule to the Markov modulated setting and  show that, under conditions, an ``averaged" $c\mu^*$ rule is asymptotically optimal, where $\mu^*$ is the ``averaged" service rate introduced in Assumption \ref{assump} (ii). Our precise  result
is given in Theorem \ref{optimal} below.  We begin with  the following definition. Let $(\p(1), \p(2), \ldots, \p(K))$ be a permutation of $(1, 2, \ldots, K)$ such that 
\be\label{ordering}
c_{\p(1)}\mu^*_{\p(1)} \geq c_{\p(2)}\mu^*_{\p(2)} \geq \cdots \geq c_{\p(K)}\mu^*_{\p(K)}. 
\ee
\begin{definition}[\textbf{$c\mu^*$ rule}]\label{optimal-policy}
 For $i\in\KK$, let
\be\label{policy1}
\mathbf{T}^n_{\p(i)}(t)=\int_0^t 1_{\{Q^n_{\p(i)}(s)>0\}}\Pi_{j=1}^{i-1}1_{\{Q^n_{\p(j)}(s)=0\}}ds.
\ee
\end{definition}
Note that $\mathbf{T}^n$ defines a priority policy in which the server always gives service priority to the nonempty queue with the largest value of $c_i\mu^*_i, i\in\KK$.

Theorem \ref{optimal} proves the asymptotic optimality of $\{\T^n\}$. Furthermore, in Theorem \ref{optimal-limit}, we will  characterize the limit of $\hatj^n(\T^n)$ in terms of the solution of a Brownian control problem which is introduced in Section \ref{bcp}. Recall that the background environment Markov process modulating the network is $\{Y^n(t) = \tilde{Y}^n(n^\nu t): t\geq 0\}.$ Consider the following three regimes for  $\nu$ and $\alpha$. 
\be\label{scaling-cases}\ba
\mbox{\bf Case 1.}  \ \ &\\
(a) \;\; &\nu>\frac{1}{2}, \;  \alpha=\frac{1}{2}. \\
(b)\;\; &0 < \nu \le 1/2, \; \alpha=\frac{1}{2}, \; \mbox{and} \; \mu^n(y)= \mu^{n,*}, \mbox{ for all } y \in \LL. \\
\mbox{\bf Case 2.} \;\;  &\nu= 0, \; \alpha=\frac{1}{2}, \; \mbox{and} \; \mu^n(y)= \mu^{n,*}, \mbox{ for all } y \in \LL. \\
\mbox{\bf Case 3.} \;\; &-1<\nu<0, \; \alpha=\frac{1-\nu}{2},  \; \mbox{and} \;  \mu^n(y)= \mu^{n,*}, \mbox{ for all } y \in \LL.
\ea\ee
Case 1 considers a situation where the jump rates of $Y^n$ are much higher than the  interarrival and service rates. For Cases 1(b),  2 and 3, the service rate $\mu^n$ is independent of $Y^n$. In Case 2, $Y^n$ has jump rates of the same order as the arrival and service rates, while in Case 3, $Y^n$ is changing slowly in comparison with the arrival and service processes. 

Throughout this work Assumptions \ref{uniferg}, \ref{assump} and \ref{ht-assump} will be taken to hold and will not be noted explicitly in statements of results.

\begin{theorem}\label{optimal}
Let $\{\mathbf{T}^n\}$ be as in Definition \ref{optimal-policy}. Then for each of the three cases in \eqref{scaling-cases}, 
\bes
\lim_{n\go\infty} \hat{J}^n(\mathbf{T}^n) = \inf\liminf_{n\go\infty}\hat{J}^n(T^n),
\ees
where the infimum is taken over all admissible control policies $\{T^n\}$.
\end{theorem}
\begin{remark}\label{remarka}
	Note that in Cases 1(b),2 and 3 the proposed control policy reduces to a standard $c\mu$-rule.  Thus Theorem \ref{optimal} shows that in these regimes, the classical $c\mu$ rule to be optimal under conditions that in particular require 
	an `averaged' form of heavy traffic condition (Assumption \ref{ht-assump}) to hold.
	\end{remark}

\begin{remark}\label{remark}
 One may conjecture that Theorem \ref{optimal} is true even when the condition that $\mu^n$ is constant is dropped  from Cases 1(b), 2 and 3. The key difficulty in treating these cases is in the proof of results analogous to Lemma \ref{ergo}. Although we are unable to treat these
cases, simulation results below suggest that the $c\mu^*$ rule performs better than the  dynamic $c\mu$ rule in these settings. 
\end{remark}
\begin{example} \label{ex}
Consider the special case of the model in Section \ref{netmod}  with, $\QQ^n = \QQ$. 
$$\KK=\{1,2\}, \;\; \LL=\{1,2\}, \; \; \pi=\left(\frac{1}{3}, \frac{2}{3}\right), \;\; c=\left(20, 25\right).$$ 
For the $n^{th}$ system,  let for $y\in\LL$,
\bes
\lambda^n_1(y)= 1 + \frac{3y}{5\sqrt{n}},\;\; \mu^n_1(y) = \frac{5}{2} + \frac{3y}{\sqrt{n}},\; \;\lambda^n_2(y)= \frac{3}{2} + \frac{3y}{5\sqrt{n}},\; \; \mu^n_2(y) = \frac{3y}{2} + \frac{3y}{\sqrt{n}}.
\ees
Then
\bes
\lambda^{n,*}_1=1 + \frac{1}{\sqrt{n}},\;\; \mu^{n,*}_1 = \frac{5}{2} + \frac{5}{\sqrt{n}},\;\; \lambda^{n,*}_2=\frac{3}{2} + \frac{1}{\sqrt{n}},\;\; \mu_2^{n,*} = \frac{5}{2} + \frac{5}{\sqrt{n}},
\ees
and
\bes
\lambda^{*}_1=1,\;\; \mu^{*}_1 = \frac{5}{2},\;\; \lambda^{*}_2=\frac{3}{2}, \;\; \mu_2^{*} = \frac{5}{2}.
\ees
Clearly Assumption \ref{ht-assump} is satisfied with $\alpha = 1/2$ and $b = (-2/5, -4/5)'$.
%
%
Note  that for all $n\in\NN$,
$
c_1\mu^{n,*}_1 < c_2\mu^{n,*}_2$ and so buffer 2 gets service priority under the $c\mu^*$ rule.
However, for the dynamic $c\mu$ rule, the priority changes according to the state of the modulating Markov process, since 
 for $n\geq 2$,
\be
 c_1\mu^n_1(1) > c_2\mu^n_2(1), \;\;  c_1\mu^n_1(2) < c_2\mu^n_2(2).
\ee
More precisely, in the dynamic $c\mu$ rule, when $Y^n(t)=1$, customers of Class $1$ have priority to be served, while when $Y^n(t)=2$, priority should be given to Class $2$ customers. The numerical study given below compares the performance of $c\mu^*$ rule with the dynamic $c\mu$ rule. The following table summarizes the simulation results for different scalings. 
We consider the cost criterion \eqref{eq:eq1050} with $\gamma = 2$.  Let $T=100$. We simulate $10$ sample paths of $\hatq^n(t)$, which are denoted as $(q^n_{1,s}(t), q^n_{2,s}(t)), s=1,2\ldots,10$, for each scaling case at discrete times $t=0.1, 0.2,\ldots, T$. We then calculate the discounted cost as follows:

$$\mbox{Cost}\; =  \frac{c_1}{10} \sum_{s=1}^{10} \sum_{t=0.1}^T e^{-\gamma t} q^n_{1,s}(t) + \frac{c_2}{10}\sum_{s=1}^{10} \sum_{t=0.1}^Te^{-\gamma t} q^n_{2,s}(t).$$
 We find in our numerical results that the discounted cost is always smaller under the $c\mu^*$ rule. 

\begin{table}[htdp]
\begin{center}
\begin{tabular}{|c|c|c|c|c|}
\hline
$n$ & $\nu$ & $\alpha$  & 
Cost ($c\mu^*$ rule) & Cost (dynamic $c\mu$ rule) \\
\hline
$25$ & $1$& $1/2$  & 
$52.70$&$55.97$ \\
\hline
$100$ & $1$& $1/2$  &
$60.18$&$61.87$\\
\hline
$25$ & $2/3$& $1/2$  & 
$72.66$&$78.57$ \\
\hline
$100$ & $2/3$& $1/2$  &
$75.07$&$77.13$\\
\hline
$25$ & $1/3$& $1/2$  & 
$63.93$&$65.04$ \\
\hline
$100$ & $1/3$& $1/2$  &
$81.60$&$84.33$\\
\hline
$25$ & $0$& $1/2$  & 
$64.91$&$69.69$\\
\hline
$100$ & $0$& $1/2$  &
$68.22$&$75.52$\\
\hline
$25$ & $-1/3$& $1/2$  &
$79.58$&$81.95$\\
\hline
$100$ & $-1/3$& $1/2$  &
$55.63$&$58.06$\\
\hline
$25$ & $-1/3$& $2/3$  &
$41.25$&$44.20$\\
\hline
$100$ & $-1/3$& $2/3$ & 
$25.80$&$27.23$\\
\hline
\end{tabular}
\end{center}
\end{table}%
Also note that only the first four rows of the table correspond to Case 1(a).  The remaining rows, although consistently show that $c\mu^*$ rule outperforms the dynamic $c\mu$ rule, correspond to regimes not covered by our 
results.
We also plot the cost at each discrete time for the average sample path (of the 10 sample paths) for a range of values of $\nu$ when $n=100$. The figures can be found in Appendix \ref{simulation}. 

\end{example}

\section{Brownian control problem}\label{bcp}

Our control policy $\{\T^n\}$ is motivated by the solution of a Brownian control problem (BCP). In this section we introduce this BCP. Roughly speaking, the BCP is obtained by taking a formal limit, as $n\to \infty$, in the equation governing the queue length evolution for the $n^{th}$ system. Recall the evolution equations for $\hat Q^n$ from \eqref{queue-hat1}.
The last term on the right side of \eqref{queue-hat1} can be rewritten as
\bes\ba
&\frac{1}{n^\alpha}\int_0^{nt}\left(\lambda^n_i(Y^n_s)-\mu^n_i(Y^n_s)\dot{T}^n_i(s)\right)ds \\
&= \mu^{n,*}_in^{1-\alpha}\left(\frac{\lambda_i^{n,*}}{\mu^{n,*}_i}-\frac{\lambda_i^*}{\mu^*_i}\right)t + \mu^{n,*}_in^{1-\alpha}\left(\frac{\lambda_i^*t}{\mu^*_i}-\bart^n_i(t)\right) \\
&\quad+ \frac{1}{n^\alpha}\int_0^{nt}\left(\lambda^n_i(Y^n_s)-\lambda^{n,*}_i\right)ds - \frac{1}{n^\alpha}\int_0^{nt}\left(\mu^n_i(Y^n_s)-\mu^{n,*}_i\right)\dott^n_i(s)ds.
\ea\ees
Let
\be\ba\label{eqn-X}
\hatx^n_i(t) &= \hatn^n_{a,i}\left(\frac{1}{n}\int_0^{nt}\lambda^n_i(Y^n_s)ds\right)-\hatn^n_{s,i}\left(\frac{1}{n}\int_0^{nt}\mu^n_i(Y^n_s)\dott^n_i(s)ds\right) \\
&\quad+ \frac{1}{n^\alpha}\int_0^{nt}\left(\lambda^n_i(Y^n_s)-\lambda^{n,*}_i\right)ds - \frac{1}{n^\alpha}\int_0^{nt}\left(\mu^n_i(Y^n_s)-\mu^{n,*}_i\right)\dott^n_i(s)ds\\
&\quad + \mu^{n,*}_in^{1-\alpha}\left(\frac{\lambda_i^{n,*}}{\mu^{n,*}_i}-\frac{\lambda_i^*}{\mu^*_i}\right)t,
\ea\ee
and
\be\label{eqn-eta}
\hat{\eta}^n_i(t) = n^{1-\alpha}\left(\frac{\lambda_i^*t}{\mu^*_i}-\bart^n_i(t)\right).
\ee
Then for $i\in\KK$ and $t\geq 0$,
\be\label{eqn-Q}
\hatq^n_i(t) = \hatx^n_i(t) + \mu_i^{n,*}\hat{\eta}^n_i(t).
\ee
We now formally take the limit as $n\to \infty$.  Write $\hatx^n=(\hatx^n_1,\ldots,\hatx^n_K)', \hat\eta^n=(\hat\eta^n_1,\ldots,\hat\eta^n_K)'$.
In Lemma \ref{conv} it is shown that, under conditions, as $n\to \infty$,
\be\label{eqn5}
\hatx^n\Go X,
\ee
where $X$ is a $K$-dimensional Brownian motion with drift $\vartheta$ and covariance matrix $\Sigma$, where $\vartheta$ and $\Sigma$ are given explicitly in Lemma \ref{conv}. 
In particular the value of $\Sigma$ is different in the three cases considered in \eqref{scaling-cases}.  Also, as $n\to \infty$, $\mu^{n,*} \to \mu^*$. 
Although in general $\hat \eta^n$ will not converge, formally taking limit in \eqref{eqn-Q}, as $n \to \infty$,  we arrive at the following BCP.
\begin{definition}[Brownian control problem (BCP)] \label{bcp-def} 
Let $X$ be a $K$-dimensional Brownian motion with drift $\vartheta$ and covariance matrix $\Sigma$ 	given on some filtered probability space
$(\bar \Omega, \bar \clf, \bar \PP, \{\bar \clf_t\})$.
The BCP is to find an $\RR^K$-valued RCLL $\{\bar \clf_t\}$ adapted stochastic process $\tileta = (\tileta_1, \tileta_2, \ldots, \tileta_K)'$ which minimizes
\bes
\EE\left(\int_0^\infty e^{-\gamma t} (c\cdot \tilq(t))dt\right)
\ees
subject to
\be\label{eqn8}\ba
& \tilq_i(t) = X_i(t) +\mu_i^*\tileta_i(t)\geq 0, \; \mbox{for all} \; i\in \KK, t\geq 0,\\
& \tili(\cdot):=\sum_{i=1}^K\tileta_i(\cdot) \; \mbox{is nondecreasing and} \; \tili(0)=0.
\ea\ee
Denote by $\tilde \cla$ the collection of all $\RR^K$-valued RCLL $\{\bar \clf_t\}$ adapted stochastic process $\tileta$ that satisfy \eqref{eqn8}.
\end{definition}
Brownian control problems of the above form were first formulated  by Harrison in \cite{harrison88}, and since then, many authors have used such control problems in the study of optimal scheduling for multiclass queuing networks in heavy traffic.
The BCP in Definition \ref{bcp-def}  has an explicit solution  given in Lemma \ref{bcp-soln} which can be proved by standard methods (see \cite{bellwill01} and references therein).
For completeness we provide the proof of Lemma \ref{bcp-soln} in Appendix \ref{bcp2}.
\begin{lemma}\label{bcp-soln}
Define for $t\geq 0$,
\be\label{opt-workload}
I^*(t) = -\inf_{0\leq s\leq t} \sum_{i=1}^K\frac{X_i(s)}{\mu^*_i},\;
W^*(t) = \sum_{i=1}^K\frac{X_i(t)}{\mu^*_i} + I^*(t),
\ee
and
\bes
\eta^*_{\p(K)}(t) = W^*(t) - \frac{X_{\p(K)}(t)}{\mu^*_{\p(K)}}, \;  \; \eta^*_j(t) = -\frac{X_j(t)}{\mu^*_j}, \; j=\p(1), \ldots, \p(K-1).
\ees
Then 
 $\eta^*$ is an optimal solution of the BCP. 
\end{lemma}
Define for $i\in\KK$, $Q^*_i$ by \eqref{eqn8} with $\tileta_i$ replaced by $\eta^*_i$. Then we have
\be\label{optimal-q}
Q^*_{\p(K)}(t) = \mu^*_{\p(K)}W^*(t)\geq 0, \; \mbox{and} \; Q^*_j(t)=0, j=\p(1),\ldots, \p(K-1).
\ee 
The optimal solution \eqref{optimal-q} suggests that the jobs should be always kept in the buffer with smallest value of $c_i\mu^*_i, i\in\KK$. This motivates the control policy defined in Definition \ref{optimal-policy}. Define
\bes
J^*= \EE\left(\int_0^\infty e^{-\gamma t} (c\cdot Q^*(t))dt\right).
\ees
From Lemma \ref{bcp-soln} we have $J^* = \inf\EE\left(\int_0^\infty e^{-\gamma t} (c\cdot \tilq(t))dt\right)$, where the infimum is taken 
over  all  $\tileta \in \tilde \cla$. In fact, we have the following result, which is proved in the next section. 
\begin{theorem}\label{optimal-limit}
	Let $\{\T^n\}_{n\ge 1}$ be as in Definition \ref{optimal-policy}.  Then for each of the three cases in \eqref{scaling-cases},
$\lim_{n\go\infty} \hatj^n(\T^n) = J^*$. 
\end{theorem}

\section{Proofs of Theorems \ref{optimal} and \ref{optimal-limit}} \label{proof}
In this section we prove Theorems \ref{optimal} and \ref{optimal-limit} . We begin with the following two lemmas, which play a crucial role in the proofs. The first lemma is a functional central limit theorem (FCLT)  for the sequence of Markov processes $\tily^n$. We first introduce some notation. 
Denote by $\clq^n$ the infinitesimal generator of $\tily^n$, namely $\clq^n: \bm(\LL) \to \bm(\LL)$ defined as
 \be\label{q-generator}
\clq^n f(y) = \sum_{j\in\LL} \qq^n_{yj}f(j), \;\; f \in \bm(\LL)\;\;\; y\in\LL.
\ee
As in Assumption \ref{uniferg}, denote by $P^{n,t}$ the transition kernel of $\tily^n$, namely, for $y\in\LL, A\subset\LL, P^{n,t}(y, A)=\PP^n(\tily^n(t)\in A \mid \tily^n(0) =y)$.
For notational simplicity we denote this probability by $\PP^n_y(\tily^n(t)\in A)$ and  by $\EE_y$ the corresponding expectation operator.  In
particular for $f\in \bm(\LL), P^{n,t}(y, f) = \EE_y(f(\tily^n_t))$. Define for $t\geq 0$ and $y\in\LL$, 
\be\label{poisson-soln1}
\lambda^{n,c}_t(y) := \EE_y\left(\int_0^t \lambda^n(\tily^n_s)ds - t\pi^n(\lambda^n)\right).
\ee
\begin{lemma}\label{fclt-Y}
	For $y\in\LL$,  $\lim_{n\go\infty}\lim_{t\go\infty}\lambda^{n,c}_t(y)$ exists.  
	Denote the limit by $\hat\lambda(y)$. Let $\{\parvs(n)\}$ be a nonnegative sequence such that $\parvs(n) \to \infty$ as $n \to \infty$.
Define for $t\geq 0$ and $n \ge 1$
\bes
G^{n}(t) = \frac{1}{\sqrt{\parvs(n)}}\int_0^{\parvs(n) t} (\lambda^n(\tily^n_s)-\lambda^{n,*})ds.
\ees
  Then $G^{n}$ converges weakly to a $K$-dimensional Brownian motion with drift $0$ and covariance matrix $\Lambda$, where for $i,j\in\KK$, 
$$\Lambda_{ij}= \int_\LL\left[(\lambda_i(y)-\lambda^*_i)\hat\lambda_j(y)+(\lambda_j(y)-\lambda^*_j)\hat\lambda_i(y)\right]\pi(dy).$$
\end{lemma}
{\bf Proof:} 
From Assumption \ref{uniferg}, there exists $a_1\in (0,\infty)$ such that for $0\leq t_1\leq t_2<\infty$, 
\bes
|\lambda^{n,c}_{t_1}(y)-\lambda^{n,c}_{t_2}(y)| \leq \int_{t_1}^{t_2} \left|P^{n,s}(y, \lambda^n)-\pi^n(\lambda^n)\right|ds \leq a_1 R \int_{t_1}^{t_2} \rho^s ds =  a_1 R \log(\rho)(\rho^{t_2}-\rho^{t_1}).
\ees
This shows that $\lim_{t\go\infty}\lambda^{n,c}_t(y)$ exists, and we denote this limit by $\hat\lambda^n(y)$. 
Note that $\hat \lambda^n(y) =  \int_0^\infty \EE_y\left(\lambda^n(\tily^n_s) - \pi^n(\lambda^n)\right)ds$.
Next note that for $y\in\LL$,
\be\label{eq:eq1201}\ba
\clq^n\hat\lambda^n(y)&= \lim_{t\go\infty}\clq^n \lambda^{n,c}_t(y) = \lim_{t\go\infty} \left(\int_0^t \sum_{j\in\LL}\mathbb{Q}_{yj}^nP^{n,s}(j,\lambda^n)ds - t\pi^n(\lambda^n)\sum_{j\in\LL}\mathbb{Q}_{yj}^n\right)\\
& =  \lim_{t\go\infty} \int_0^t \sum_{j\in\LL}\mathbb{Q}_{yj}^nP^{n,s}(j,\lambda^n)ds
= \lim_{t\go\infty}\int_0^t \frac{dP^{n,s}}{ds} (y,\lambda^n)ds
 = \pi^n(\lambda^n) - \lambda^n(y).
\ea\ee
Noting that $\mathbb{Q}^n$ converges to $\mathbb{Q}$, it follows that (cf. \cite[Theorem 4.2.5]{ethier-kurtz86}), if $\tily^n (0)$ has limit distribution $\tilde\nu$ on $\LL$, then there exists an $\LL$-valued Markov process $\tily$ with initial distribution $\tilde\nu$ and rate matrix $\mathbb{Q}$ such that $\tily^n\Go\tily$. Therefore, by 
Assumption \ref{assump} (i), Assumption \ref{uniferg} and  
dominated convergence theorem, as $n\go\infty$,
\bes
\hat\lambda^n(y) = \int_0^\infty \EE_y\left(\lambda^n(\tily^n_s) - \pi^n(\lambda^n)\right)ds \go \int_0^\infty\EE_y\left( \lambda(\tily_s) - \pi(\lambda)\right)ds,
\ees
This proves the first statement in the lemma, in fact
\bes
\hat\lambda(y) = \int_0^\infty \EE_y\left(\lambda(\tily_s) - \pi(\lambda)\right)ds.
\ees
Following the proof of
\eqref{eq:eq1201} one can
check that $\clq\hat\lambda(y) = \pi(\lambda) - \lambda(y)$ for $y\in\LL$, where $\clq$ is defined as in \eqref{q-generator} with $\QQ^n$ replaced by $\QQ$. Next observe that for each $n$,
\be\label{martingale}
M^n_t :=\hat\lambda^n(\tily^n_t) - \hat\lambda^n(\tily^n_0) -\int_0^t \clq^n \hat\lambda^n(\tily^n_s) ds, \;\; t\geq 0, 
\ee
is a $\clg^n_{t} = \sigma \{\tily^n_s: s \le t\}$ square integrable martingale.  Furthermore
\be\label{difference}
\sup_{t\geq 0}|G^{n}(t) - \parvs(n)^{-\frac{1}{2}}M^n_{\parvs(n)t}| \leq \parvs(n)^{-\frac{1}{2}} \sup_{t\geq 0}|\hat\lambda^n(\tily^n(\parvs(n) t)) - \hat\lambda^n(\tily^n(0)) | \go 0. 
\ee
Thus it suffices to show 
\be\label{fclt-m}
\parvs(n)^{-\frac{1}{2}}M^n_{\parvs(n)\cdot} \mbox{ {\small converges weakly to a   Brownian motion with drift $0$ and covariance matrix $\Lambda$.}}
\ee
Proof of \eqref{fclt-m} follows from classical functional central limit theorems for martingales. For completeness, we provide the  details  in   Appendix \ref{lemma}. \ink

\begin{lemma}\label{ergo}
Let $\{f^n\}_{n\ge 1}$, $f^n: \LL\go \RR$, be a sequence of functions such that $\sup_{n\in\NN, y\in\LL}|f^n(y)|<\infty$. Define for $t\geq 0$ and $i\in\KK$,
\bes
\Phi^{f^n}(t) = \frac{1}{\sqrt{n}}\int_0^{nt} (f^n(Y^n_s)-f^{n,*})\dott^n_i(s)ds,
\ees
where $f^{n,*}= \pi^n(f^n)$. Suppose that $\nu > 1/2$.
Then, there exists $\kappa_1 \in (0, \infty)$ such that, for all $t \ge 0$
\be\label{suffice} \EE\left[ \sup_{0 \le s \le t} |\Phi^{f^n}(s)| \right] \le \frac{\kappa_1 t}{n^{\nu -1/2}}.\ee
In particular, as $n\go\infty$, $\Phi^{f^n} \Go 0.$
\end{lemma}
{\bf Proof:} 
Note that for $t\geq 0$,
\be\label{first-eqn}
\ba
\frac{1}{\sqrt{n}}\int_0^{nt}(f^n(Y^n_s)-f^{n,*})\dot{T}^n_i(s)ds &= \frac{1}{n^{\nu+\frac{1}{2}}}\int_0^{n^{\nu+1}t} (f^n(\tily^n_s)-f^{n,*})\dot{T}^n_i(s/{n^\nu})ds\\
&=\frac{1}{n^{\nu+\frac{1}{2}}}\sum_{k=0}^{\lfloor n^{\nu+1} t\rfloor-1}\int_k^{k+1} (f^n(\tily^n_s)-f^{n,*})\dot{T}^n_i(s/{n^\nu})ds \\
&\quad+ \frac{1}{n^{\nu+\frac{1}{2}}}\int_{\lfloor n^{\nu+1} t\rfloor}^{n^{\nu+1} t} (f^n(\tily^n_s)-f^{n,*})\dot{T}^n_i(s/n^\nu)ds.
\ea
\ee
For the last  expression in \eqref{first-eqn}, note that there exists $a_1\in(0,\infty)$ such that for all $t\geq 0$,
\be\label{mainlemma-ineqn1}
\left|\frac{1}{n^{\nu+\frac{1}{2}}}\int_{\lfloor n^{\nu+1} t\rfloor}^{n^{\nu+1} t} (f^n(\tily^n_s)-f^{n,*})\dot{T}^n_i(s/{n^\nu})ds\right|\leq \frac{a_1}{n^{\nu+\frac{1}{2}}}.
\ee
In what follows, we set $\dott^n_i(t)=0$ when $t < 0$. Then 
\bes\ba
&\frac{1}{n^{\nu+\frac{1}{2}}}\sum_{k=0}^{\lfloor n^{\nu+1} t\rfloor-1}\int_k^{k+1} (f^n(\tily^n_s)-f^{n,*})\dot{T}^n_i(s/{n^\nu})ds \\
 &= \frac{1}{n^{\nu+\frac{1}{2}}}\sum_{k=0}^{\lfloor n^{\nu+1} t\rfloor-1}\dot{T}^n_i\left(\frac{k-2}{n^\nu}\right)\int_k^{k+1} (f^n(\tily^n_s)-f^{n,*})ds\\
&\quad+\frac{1}{n^{\nu+\frac{1}{2}}}\sum_{k=0}^{\lfloor n^{\nu+1} t\rfloor-1}\int_k^{k+1} (f^n(\tily^n_s)-f^{n,*})\left(\dot{T}^n_i\left(\frac{s}{n^\nu}\right)-\dot{T}^n_i\left(\frac{k-2}{n^\nu}\right)\right)ds.
\ea\ees
Next note that, there exists $a_2 \in (0, \infty)$ such that for $2\leq k\leq \lfloor n^{\nu+1} t\rfloor-1$,
\be\label{Tdot-est}\ba
&\EE\left(\sup_{k\leq s\leq k+1}\left|\dot{T}^n_i\left(\frac{s}{n^\nu}\right)-\dot{T}^n_i\left(\frac{k-2}{n^\nu}\right)\right|\right)\\
&\leq \PP\left(\sup_{k-2\leq s\leq k+1}\left|\dot{T}^n_i\left(\frac{s}{n^\nu}\right)-\dot{T}^n_i\left(\frac{k-2}{n^\nu}\right)\right|>0\right)  \le \frac{a_2}{n^{\nu}},
\ea\ee
where the last inequality is a consequence of the fact that
 the arrival and service rates are $\clo(1)$, and by Condition \ref{policy-cond} (i), $\dott^n$ does not change values between two successive arrivals and/or service completions. 
Thus we can find $a_3\in(0,\infty)$ such that for all $t\geq 0$,
\be\label{mainlemma-ineqn2}
\EE\left(\sup_{0\leq s\leq t}\left|\frac{1}{n^{\nu+\frac{1}{2}}}\sum_{k=0}^{\lfloor n^{\nu+1} s\rfloor-1}\int_k^{k+1} (f^n(\tily^n_u)-f^{n,*})\left(\dot{T}^n_i\left(\frac{u}{n^\nu}\right)-\dot{T}^n_i\left(\frac{k-2}{n^\nu}\right)\right)du\right|\right)\leq \frac{a_3t}{n^{\nu-\frac{1}{2}}}.
\ee
Define $g^n \in \bm(\LL)$ as
$$g^n(y) = \int_0^{\infty} \EE_y\left (f^n(\tily^n(s))) - \pi^n(f^n)\right) ds.$$
From Assumption \ref{uniferg} and calculations as in the proof of Lemma \ref{fclt-Y}  it follows that
$\max_{y\in \LL}\sup_n |g^n(y)| < \infty$.  Also,
%
%
%
%
\bes
\tilde M^n(s):= g^n(\tilde{Y}^n_s)-g^n(\tilde{Y}^n_0)+\int_0^s (f^n(\tilde{Y}^n_u)-f^{n,*})du, \; s \ge 0
\ees
is a $\{\clg^n_t\}$ martingale, where $\clg_t^n = \sigma \{\tily^n_s: s \le t\}$. For $0\leq k\leq \lfloor n^{\nu+1} t\rfloor-1$, let $\xi_k^n = \tilde M^n(k+1)-\tilde M^n(k).$
Then
\begin{align*} 
&\frac{-1}{n^{\nu+\frac{1}{2}}}\sum_{k=0}^{\lfloor n^{\nu+1} t\rfloor-1}\dot{T}^n_i\left(\frac{k-2}{n^\nu}\right)\int_k^{k+1} (f^n(\tily^n_s)-f^{n,*})ds \\
& = \frac{1}{n^{\nu+\frac{1}{2}}}\sum_{k=0}^{\lfloor n^{\nu+1} t\rfloor-1}\dot{T}^n_i\left(\frac{k-2}{n^\nu}\right)\left[g^n(\tilyn_{k+1})-g^n(\tilyn_k)-\xi^n_k\right]\\
\end{align*}
The last expression equals
\be \label{eq:eq1308}\ba
& \frac{1}{n^{\nu+\frac{1}{2}}}\sum_{k=0}^{\lfloor n^{\nu+1} t\rfloor-1}\left[\dot{T}^n_i\left(\frac{k-2}{n^\nu}\right)g^n(\tilyn_{k+1})- \dot{T}^n_i\left(\frac{k-3}{n^\nu}\right)g^n(\tilyn_k)\right] \\
&+ \frac{1}{n^{\nu+\frac{1}{2}}}\sum_{k=0}^{\lfloor n^{\nu+1} t\rfloor-1}\left[\dot{T}^n_i\left(\frac{k-3}{n^\nu}\right)-\dot{T}^n_i\left(\frac{k-2}{n^\nu}\right)\right]g^n(\tilyn_k)- \frac{1}{n^{\nu+\frac{1}{2}}}\sum_{k=0}^{\lfloor n^{\nu+1} t\rfloor-1}\dot{T}^n_i\left(\frac{k-2}{n^\nu}\right)\xi^n_k.
\ea\ee
The first term in the above display is a telescoping sum and since $\{g^n\}$ is uniformly bounded, we  can find $a_4\in(0,\infty)$ such that for all $t\geq 0$,
\be\label{mainlemma-ineqn3}
\left|\frac{1}{n^{\nu+\frac{1}{2}}}\sum_{k=0}^{\lfloor n^{\nu+1} t\rfloor-1}\left[\dot{T}^n_i\left(\frac{k-2}{n^\nu}\right)g^n(\tilyn_{k+1})- \dot{T}^n_i\left(\frac{k-3}{n^\nu}\right)g^n(\tilyn_k)\right] \right| \leq \frac{a_4}{n^{\nu+\frac{1}{2}}}.
\ee
Also, from \eqref{Tdot-est}, there exists $a_5\in(0,\infty)$ such that for all $t\geq 0$,
\be\label{mainlemma-ineqn4}
\EE\left(\sup_{0\leq s\leq t}\left| \frac{1}{n^{\nu+\frac{1}{2}}}\sum_{k=0}^{\lfloor n^{\nu+1} s\rfloor-1}\left[\dot{T}^n_i\left(\frac{k-3}{n^\nu}\right)-\dot{T}^n_i\left(\frac{k-2}{n^\nu}\right)\right]g^n(\tilyn_k)\right|\right)\leq \frac{a_5t}{n^{\nu-\frac{1}{2}}}.
\ee
Now we consider the last term in \eqref{eq:eq1308}.

For $0\leq k\leq \lfloor n^{\nu+1} t\rfloor-1$, let $\tilxi^n_k= \dot{T}^n_i\left(\frac{k-2}{n^\nu}\right)\xi^n_k$. By Condition \ref{policy-cond} (ii), 
$\EE\left(\tilxi^n_k |\tilde{\clf}^n_k\right)=0$ and so $\tilxi^n_k$ is a martingale difference. Furthermore, $\tilxi^n_k$ is uniformly bounded in $n$ and $k$. Applying Doob's inequality, there exists $a_6\in(0,\infty)$ such that for $t\geq 0$,
\be\label{mainlemma-ineqn5}\ba
&\EE\left(\sup_{0\leq s\leq t}\left|\frac{1}{n^{\nu+\frac{1}{2}}}\sum_{k=0}^{\lfloor n^{\nu+1} s\rfloor-1}\dot{T}^n_i\left(\frac{k-2}{n^\nu}\right)\xi^n_k\right|\right)^2\leq \frac{4}{n^{2\nu+1}}\EE\left(\sum_{k=0}^{\lfloor n^{\nu+1} t\rfloor-1}\dot{T}^n_i\left(\frac{k-2}{n^\nu}\right)\xi^n_k\right)^2\\
&\quad=\frac{4}{n^{2\nu+1}}\sum_{k=0}^{\lfloor n^{\nu+1} t\rfloor-1}\EE\left(\dot{T}^n_i\left(\frac{k-2}{n^\nu}\right)\xi^n_k\right)^2 \leq \frac{a_6t}{n^{\nu}}.
\ea\ee
The result follows on combining \eqref{mainlemma-ineqn1}, \eqref{mainlemma-ineqn2}, and \eqref{mainlemma-ineqn3} -- \eqref{mainlemma-ineqn5}. \ink

The following lemma gives a functional law of large numbers (FLLN) and a 
 FCLT.
For a sequence of admissible policies $\{T^n\}$, $i\in\KK$ and $t\geq 0$, 
let
\be\label{job-completion-proc}
D^n_i(t) = N^n_{s,i}\left(\int_0^t\mu^n_i(Y^n_s)\dott^n_i(s)ds\right).
\ee
Recall $\{A^n_i\}$ from equation \eqref{eq:eq5.1}.
Write $A^n=(A^n_1,\ldots,A^n_K)'$ and $D^n=(D^n_1,\ldots, D_K)'$. Define the fluid scaled processes $\bara^n(t) = n^{-1}A^n(nt)$ and $\bard^n(t) = n^{-1}D^n(nt)$. 
Let $\iota:[0,\infty)\go[0,\infty)$ be the identity map, and define for $t\geq 0$,
\be\label{T-star}
\bart^*(t) = \left(\frac{\lambda_1^*}{\mu^*_1}, \frac{\lambda_2^*}{\mu^*_2}, \ldots, \frac{\lambda_K^*}{\mu^*_K}\right)' t.
\ee
Finally, recall  $\hatx^{n}$ defined in \eqref{eqn-X}.

\begin{lemma}\label{conv} 
Let $\{T^n\}$ be a sequence of admissible scheduling policies such that
\be \label{eq:eq1556}
\underline{J}(\{T^n\})= \liminf_{n\rightarrow\infty}\hat{J}^n(T^n) <\infty.
\ee
Let $\{n'\}$ be a subsequence of $\{n\}$ along which the above $\liminf$ is attained. 
\bi
\item[\rm (i)] The following FLLN  holds: As $n' \to \infty$
\bes
(\barq^{n'}, \bara^{n'}, \bard^{n'}, \bart^{n'}, \bari^{n'}) \Rightarrow (0, \lambda^*\iota, \lambda^*\iota, \bart^*, 0).
\ees
\item[\rm (ii)] The following FCLT holds:  As $n' \to \infty$
\bes
\hatx^{n'} \Go X, 
\ees
where $X$ is a Brownian motion with drift $\vartheta= (\mu_1^*b_1,\ldots,\mu_K^*b_K)'$ and covariance matrix $\Sigma$
which is given in the three regimes as follows: 
\bes\ba
\mbox{\bf Case 1.}  \;\; & \Sigma = \diag(2\lambda^*). \\
\mbox{\bf Case 2.} \;\;  & \Sigma = \diag(2\lambda^*) + \Lambda. \\
\mbox{\bf Case 3.} \;\; & \Sigma = \Lambda.
\ea\ees
\ei
\end{lemma}
{\bf Proof:} 
With abuse of notation, we relabel the subsequence $\{n'\}$ once more as $\{n\}$ and with this relabeling we have that 
\be\label{finite-cost}
\lim_{n\rightarrow\infty} \hat{J}^{n}(T^{n}) = \underline{J}(\{T^{n}\})<\infty.
\ee

Part (i). Let $i\in\KK$. From functional law of large numbers for Poisson processes, we have
\be\label{poisson}
\barn_{a,i}^n(\cdot)\Go \iota(\cdot), \; \mbox{and} \; \barn_{s,i}^n(\cdot)\Go \iota(\cdot).
\ee
Next note that
\bes
 \frac{1}{n}\int_0^{nt} (\lambda^n_i(Y^n_s) - \lambda^{n,*}_i)ds = \frac{1}{n^{\nu+1}}\int_0^{tn^{\nu+1}} (\lambda^n_i(\tily^n_s) - \lambda^{n,*}_i)ds .
\ees
Thus from Lemma \ref{fclt-Y}, we have that 
\be\label{eqnA1}
 \frac{1}{n}\int_0^{n\cdot} (\lambda^n_i(Y^n_s) - \lambda^{n,*}_i)ds\Go 0,
\ee
and thus  by Assumption \ref{assump}(i),
\be\label{eqnA3}
\frac{1}{n}\int_0^{n\cdot} \lambda^n_i(Y^n_s)ds \Go \lambda^*_i \iota(\cdot).
\ee
Using the random change of time theorem (see \cite[Lemma 3.14.1]{billingsley99}), we have that
\be\label{conv-A}
\bara^n_i(\cdot)\Go \lambda_i^* \iota(\cdot).
\ee
In cases 1(b), 2 and 3 (i.e. when $\mu^n(Y^n)\equiv\mu^{n,*}$), by Assumption \ref{assump}(ii), we have that for $t\geq 0$,
\be\label{eqnA5}
\EE\left(\sup_{0\leq s\leq t}\left|\frac{1}{n}\int_0^{ns} \mu^n_i(Y^n_u)\dott^n_i(u)du - \mu^*_i \bart^n_i(s)\right|\right)\leq |\mu^{n,*}_i-\mu^*_i|t \go 0.
\ee
Furthermore, in case 1(a) (i.e.  when $\nu>1/2$ and $\mu^n$ is modulated by $Y^n$), using Assumption \ref{assump}(ii) and Lemma \ref{ergo} (see \eqref{suffice}), we have the same convergence result as above. Combining \eqref{eqnA5} and \eqref{poisson}, we have that
\be\label{tight-D}
\bard^n_i \; \mbox{is} \; C\mbox{-tight}.
\ee
From  \eqref{lip}, 
\be\label{tight-T}
\bart^n_i \; \mbox{is} \; C\mbox{-tight}.
\ee
We next note that from \eqref{queue-fluid} and \eqref{idle},
\be\label{eqnA7}
\barq^n(t) = \bara^n(t) - \bard^n(t),
\ee
and
\be\label{eqnA8}
\bari^n(t) = t -\sum_{i=1}^K \bart^n_i(t), \; t\geq 0.
\ee
From \eqref{conv-A}, \eqref{tight-T} and \eqref{tight-D}, we now have 
\be\label{tight-ADT}
(\barq^n, \bari^n, \bara^n, \bard^n, \bart^n) \; \mbox{is} \; C\mbox{-tight}.
\ee
Let $(\barq, \bari, \bara, \bard, \bart)$ be a limit point,  along a subsequence $\{n''\}$ of $\{n\}$. By Skorohod representation theorem, we can assume that this convergence holds almost surely, uniformly on compacts. From \eqref{finite-cost} and using Fatou's lemma, we have
\be\ba\label{q-bar}
0 = \lim_{n''\rightarrow\infty} \frac{1}{\sqrt{n''}}\hat{J}^{n''}(T^{n''}) &\geq \EE\left(\int_0^\infty e^{-\gamma t} \left[\liminf_{n''\rightarrow\infty} h\cdot \barq^{n''}(t)\right]dt\right)\\
&= \EE\left(\int_0^\infty e^{-\gamma t} \left[h\cdot \barq(t)\right]dt\right).
\ea\ee
Since $h>0$ and $\barq$ has continuous sample paths, we have $\barq\equiv0$. By \eqref{conv-A} and \eqref{eqnA7},
\be\label{eq:eq1333}
\bard^{n'}(\cdot)\Go \lambda^*\iota(\cdot)
\ee
Also, from \eqref{poisson} and \eqref{eqnA5}, we have that $\bard^{n''}_i\Go \mu^*_i\bart_i$. Therefore, 
\be\label{eqnA9}
\bart_i(t) = \frac{\lambda_i^*}{\mu^*_i} t, \;\; t\geq 0.
\ee
Finally, \eqref{eqnA9} and heavy traffic condition \eqref{hr1} imply $\bari^{n''}\Go 0$. This completes Part (i).

We now consider part (ii). From \eqref{eqn-X} we have that 
\be\label{eqn-X1}\ba
\hatx^n_i(t) 
&= \hatn^n_{a,i}\left(\frac{1}{n}\int_0^{nt}\lambda^n_i(Y^n_s)ds\right)-\hatn^n_{s,i}\left(\frac{1}{n}\int_0^{nt}\mu^n_i(Y^n_s)\dott^n_i(s)ds\right) \\
&\quad+ n^{-(\alpha+\frac{\nu-1}{2})}G_i^{n^{\nu+1}}(t) - n^{-\alpha+\frac{1}{2}}\Phi^{\mu^n_i}(t) +  \mu^{n,*}_in^{1-\alpha}\left(\frac{\lambda_i^{n,*}}{\mu^{n,*}_i}-\frac{\lambda_i^*}{\mu^*_i}\right)t.
\ea\ee
We now treat the three cases separately. \\
{\bf Case 1.} From functional central limit theorem for Poisson processes, we have that 
\be\label{poisson2}
(\hatn^n_{a,i}, \hatn^n_{s,i})_{i=1}^K \Go (\hatn_{a,i},\hatn_{s,i})_{i=1}^K,
\ee 
where $\hatn_{a,i}$ and $\hatn_{s,i}$, $i=1, \cdots K$ are $1$-dimensional mutually independent standard Brownian motions. From \eqref{eqnA3}, \eqref{eqnA5}, and \eqref{eqnA9}, we have
\be\label{poisson3}
\left (\hatn^n_{a,i}\left(\frac{1}{n}\int_0^{n\cdot}\lambda^n_i(Y^n_s)ds\right)-\hatn^n_{s,i}\left(\frac{1}{n}\int_0^{n\cdot}\mu^n_i(Y^n_s)\dott^n_i(s)ds\right)
\right )_{i=1}^K\Go (W_i)_{i=1}^K,
\ee
where $W_i$, $i= 1, \cdots K$ are  $1$-dimensional mutually independent Brownian motions with drift $0$ and variance $2\lambda^*_i$. 
In case 1(a) (i.e. when $\nu>1/2, \alpha=1/2$ and $\mu^n$ is modulated by $Y^n$), by Lemmas \ref{fclt-Y} and \ref{ergo}, we have $n^{-\frac{\nu}{2}}G_i^{n^{\nu+1}} + \Phi^{\mu^n_i}\Go 0$. In case 1(b) (i.e. when $\nu>0, \alpha=1/2$ and $\mu^n\equiv \mu^{n,*}$), we have the same convergence by Lemma \ref{fclt-Y}. Thus, using \eqref{hr2}, $(\hatx^n_i)_{i=1}^K$ converges weakly to $(X_i)_{i=1}^K$, where $X_i$, $i=1, \cdots K$ are mutually independent Brownian motions with drift $\mu^*_ib_i$ and covariance $2\lambda^*_i$. \\
{\bf Case 2.} Here $\alpha=\frac{1}{2}$ and $\nu=0$. As in  Case 1, we have the convergence in  \eqref{poisson3}. 
Also note that $\Phi^{\mu^n_i}\equiv 0$ since $\mu^n_i\equiv \mu^{n,*}_i$.  Furthermore, by Lemma \ref{fclt-Y}, we have 
\be\label{fclt}
G^{n^{\nu+1}}  \Go \tilde{W},
\ee
where $\tilde{W}$ is a $K$-dimensional Brownian motion with drift $0$ and covariance matrix $\Lambda$. 
The result now follows on recalling that, $Y^n, N^n_{a,i}, N^n_{s,j}, i,j\in\KK$ are mutually independent. \\
{\bf Case 3.} Here $\alpha=\frac{1-\nu}{2}$ and $-1<\nu<0$. From Lemma \ref{fclt-Y}, the convergence in \eqref{fclt} continues to hold. Also,
 since $\alpha> \frac{1}{2}$, we have 
\bes
\hatn^n_{a,i}\left(\frac{1}{n}\int_0^{n\cdot}\lambda^n_i(Y^n_s)ds\right)-\hatn^n_{s,i}\left(\frac{1}{n}\int_0^{n\cdot}\mu^n_i(Y^n_s)\dott^n_i(s)ds\right)
\Go 0.
\ees
Finally, $\Phi^{\mu^n_i}\equiv 0$. The result follows. \ink

We  recall below the definition and basic properties of   the $1$-dimensional Skorohod map. Let $\mathcal{D}_0([0,\infty), \RR)=\{x\in \mathcal{D}([0,\infty): \RR): x(0)\geq 0\}$. 

\begin{definition}[$1$-dimensional Skorohod Problem (SP)]\label{sp-def}
Let $x\in \mathcal{D}_0([0,\infty), \RR)$. A pair $(z,y)\in\mathcal{D}([0,\infty): \RR_+)\times \mathcal{D}([0,\infty): \RR_+)$ is a solution of the Skorohod problem for $x$ if the following hold.
\bi
\item[\rm (i)] For all $t\geq 0, z(t)=x(t)+y(t)\geq 0$.
\item[\rm (ii)] $y$ satisfies the following: (a) $y(0)=0$, (b) $y(\cdot)$ is nondecreasing, and (c) $y(\cdot)$ increases only when $z(\cdot)=0$, that is, 
$\int_0^\infty z(t)dy(t) = 0.$
\ei
We write $z = \Gamma(x)$ and refer to the map $\Gamma: \mathcal{D}_0([0,\infty), \RR) \to \mathcal{D}([0,\infty): \RR_+)$ as the Skorohod map.
\end{definition}
The following proposition summarizes some well known properties of the $1$-dimensional SP.

\begin{proposition}\label{sp}\hfill
\bi
\item[\rm (i)] Let $x\in \mathcal{D}_0([0,\infty), \RR)$. Then there exists a unique solution $(z,y)\in\mathcal{D}([0,\infty): \RR_+)\times \mathcal{D}([0,\infty): \RR_+)$ of the SP for $x$, which is given as follows:
\bes
y(t) = - \inf_{0\leq s\leq t} x(s), \; \mbox{and} \; z(t) = x(t) - \inf_{0\leq s\leq t} x(s), \; t \ge 0.
\ees
\item[\rm (ii)] The  Skorohod map $\Gamma$ is Lipschitz continuous in the following sense: There exists $\kappa\in (0,\infty)$ such that for all $t\geq 0$ and $x_1, x_2\in \mathcal{D}_0([0,\infty), \RR)$, 
$$\sup_{0\leq s\leq t}|\Gamma(x_1)(s)-\Gamma(x_2)(s)|\leq \kappa \sup_{0\leq s\leq t}|x_1(s)-x_2(s)|.$$
\item[\rm (iii)] Fix $x\in \mathcal{D}_0([0,\infty), \RR)$. Let $(z,y)\in\mathcal{D}([0,\infty): \RR_+)\times \mathcal{D}([0,\infty): \RR_+)$ be such that 
\bi
\item[\rm (a)] $z(t)=x(t)+y(t)\geq 0, \; t\geq 0$,
\item[\rm (b)] $y(\cdot)$ is nondecreasing with $y(0)=0$.
\ei
Then  $z(t)\geq \Gamma(x)(t), \; t\geq 0$.
\ei
\end{proposition}
We next present two lemmas that will be used in the proofs of Theorems \ref{optimal} and \ref{optimal-limit}.   Recall the control policy $\T^n$  introduced in Definition \ref{optimal-policy}. 
 Lemma \ref{conv4} says that the long run average of $\T^n$ approaches  $\bart^*$ and  Lemma \ref{conv3} shows that under $\T^n$, all the jobs are kept in buffer $\p(K)$, i.e., the buffer with smallest value of $c_i\mu^*_i, i\in\KK$, asymptotically.


%
\begin{lemma}\label{conv4} As $n\go\infty$, $\bar{\T}^n\Go \bart^*$.
\end{lemma}
{\bf Proof:} 
In the proof, all the quantities, e.g. queue lengths, idle times etc.,  are considered under the control policy $\{\T^n\}$.

Define, for $n \ge 1$,
$
W^n(t) = \sum_{i=1}^K\frac{Q^n_i(t)}{\mu^{n,*}_i},
$
and let
\be\label{diff-wl}
\hatw^n(t) =  \frac{W^n(nt)}{n^{\alpha}} = \sum_{i=1}^K\frac{\hatq^n_i(t)}{\mu^{n,*}_i} = \sum_{i=1}^K\frac{\hatx^{n}_i(t)}{\mu_i^{n,*}} + \hati^{n}(t).
\ee
Note that 
\be\label{hati}
\hati^n(t) = \frac{1}{n^\alpha}\left(nt-\sum_{i=1}^K\T^n_i(nt)\right)= n^{1-\alpha}\int_0^t 1_{\{\hatq^n_i(s)=0,i\in\KK\}}ds = n^{1-\alpha}\int_0^t 1_{\{\hatw^n(s)=0\}}ds.
\ee
Clearly $\hati^n(0)=0$,  $\hati(\cdot)$ is nondecreasing and increases only when $\hatw^n(\cdot)=0$. Therefore we have (see Definition \ref{sp-def}),  for $t\geq 0$,
\be\label{workload1}
\hatw^{n}(t) = \Gamma\left(\sum_{i=1}^K\frac{\hatx^{n}_i}{\mu_i^{n,*}}\right)(t).
\ee
 Using the Lipschitz  property  in Proposition \ref{sp} (ii),  for $t\geq 0$,
\be\label{est8}
\sup_{0\leq s\leq t}\hatw^n(s)\leq \sum_{i=1}^K\frac{\kappa}{\mu_i^{n,*}}\sup_{0\leq s\leq t}\left|\hatx^{n}_i(s)\right|.
\ee
Fix $i\in\KK$. 
Then from \eqref{eqn-X1}, we have for all $t\geq 0$,
\bes\ba
\sup_{0\leq s\leq t}\left|\hatx^{n}_i(s)\right|&\leq \sup_{0\leq s\leq t}\left|\hatn_{a,i}^n\left(\frac{1}{n}\int_0^{ns}\lambda^n_i(Y^n_u)du\right)\right|+ \sup_{0\leq s\leq t}\left|\hatn_{s,i}^n\left(\frac{1}{n}\int_0^{ns}\mu^n_i(Y^n_u)\dot{\T}^n_i(u)du\right)\right|\\
&\quad+ \sup_{0\leq s\leq t}\left|\mu^{n,*}_in^{1-\alpha}\left(\frac{\lambda^{n,*}_i}{\mu^{n,*}_i}-\frac{\lambda^*_i}{\mu^*_i}\right) s\right|+n^{-\alpha+\frac{1}{2}}\sup_{0\leq s\leq t}\left|\Phi^{\mu^n_i}(s)\right| \\
&\quad + n^{-(\alpha + \frac{\nu-1}{2})} \sup_{0\leq s\leq t}\left|G^{n^{\nu+1}}_i(s)\right|.
\ea\ees
For any $t\geq 0$, we have 
\bes
\max\left\{\frac{1}{n}\int_0^{nt}\mu^n_i(Y^n_u)\dot{\T}^n_i(u)du, \; \frac{1}{n}\int_0^{nt}\lambda^n_i(Y^n_u)du\right\}\leq \max_{l\in\LL, n\in\NN} (\lambda^n_i(l)+\mu^n_i(l)) t.
\ees
Let $a_1:= \max_{l\in\LL, n\in\NN} (\lambda^n_i(l)+\mu^n_i(l))$. From Assumption \ref{assump} (i) and (ii) we have (without loss of generality)  that $a_1<\infty$. 
Using Doob's inequality, we have for some $a_2\in(0,\infty)$,
\be\label{est-p1}\ba
\EE\left(\sup_{0\leq s\leq t}\left|\hatn_{a,i}^n\left(\frac{1}{n}\int_0^{ns}\lambda^n_i(Y^n_u)du\right)\right|\right)^2 &\leq \EE\left(\sup_{0\leq s\leq  a_1 t}\left|\hatn_{a,i}^n\left(s\right)\right|\right)^2 \leq 4\EE\left(\hatn_{a,i}^n\left(a_1 t\right)\right)^2\\
& = \frac{4a_1n t}{n^{2\alpha}} \leq a_2 t.
\ea\ee
Similarly, for some $a_3\in(0,\infty)$, we have 
\be\label{est-p2}
\EE\left(\sup_{0\leq s\leq t}\left|\hatn_{s,i}^n\left(\frac{1}{n}\int_0^{ns}\mu^n_i(Y^n_u)\dot{\T}^n_i(u)du\right)\right|\right)^2 \leq a_3 t.
\ee
Also, by Assumption \ref{ht-assump},
\be\label{est-drift}
\sup_{0\leq s\leq t}\left|\mu^{n,*}_in^{1-\alpha}\left(\frac{\lambda^{n,*}_i}{\mu^{n,*}_i}-\frac{\lambda^*_i}{\mu^*_i}\right) s\right|\go \mu^*_i|b_i|t.
\ee
Next from  Lemma \ref{ergo} (see \eqref{suffice}), we have in case 1(a), for some $a_4\in (0,\infty)$
\be\label{est-Phi}
\EE\left(n^{-\alpha+\frac{1}{2}}\sup_{0\leq s\leq t}\left|\Phi^{\mu^n_i}(s)\right|\right)\leq \frac{a_4t}{n^{\nu-\frac{1}{2}}}. 
\ee
We note that $\Phi^{\mu^n_i}\equiv 0$ in cases 1(b), 2 and 3, and so the above estimate \eqref{est-Phi} holds for all cases.
Finally, we show that for some $a_5\in(0,\infty)$,
\be\label{est-G}
\EE\left(n^{-(\alpha + \frac{\nu-1}{2})}\sup_{0\leq s\leq t}\left|G^{n^{\nu+1}}_i(s)\right|\right)^2\leq a_5 (t+1). 
\ee
For this, from \eqref{difference} and the fact that $\alpha + \frac{\nu-1}{2} \ge 0$, it suffices to show that there exists $a_6\in(0,\infty)$ such that for all $n\in\NN$,  
\bes
\EE\left(\sup_{0\leq s\leq t}\left|n^{-\frac{\nu+1}{2}}M^n_{n^{\nu+1}s}\right|\right)^2\leq a_6 (t+1). 
\ees
Using Doob's inequality and a standard martingale argument, we have that
\bes\ba
&\EE\left(\sup_{0\leq s\leq t}\left|n^{-\frac{\nu+1}{2}}M^n_{n^{\nu+1}s}\right|\right)^2 \leq 4 n^{-(\nu+1)} \EE\left(\left|M^n_{n^{\nu+1}t}\right|^2\right) \\
& = 4n^{-(\nu+1)} \sum_{1\leq k\leq \lfloor n^{\nu+1}t\rfloor} \EE\left(\left|M^n_k\right|^2-\left|M^n_{k-1}\right|^2\right) + 4n^{-(\nu+1)}\EE(|M^n_{n^{\nu+1}t}|^2 - |M^n_{\lfloor n^{\nu+1}t\rfloor}|^2) \\
& = 4n^{-(\nu+1)} \sum_{1\leq k\leq \lfloor n^{\nu+1}t \rfloor } \EE\left(\EE\left(\left.\left|M^n_k-M^n_{k-1}\right|^2\right| \clg^n_{k-1}\right)\right) \\
&\quad+ 4n^{-(\nu+1)}\EE\left(\EE\left(\left.|M^n_{n^{\nu+1}t}- M^n_{\lfloor n^{\nu+1}t\rfloor}|^2\right|\clg^n_{\lfloor n^{\nu+1}t\rfloor}\right)\right)\\
& \leq a_7 (t+1),
\ea\ees
where the last estimate uses the fact that $|M^n_k-M^n_{k-1}|$ is uniformly bounded in $n$ and $k$. Combining \eqref{est-p1} -- \eqref{est-G}, we have for some $a_8\in(0,\infty)$,
\be\label{est9}
\EE\left(\sup_{0\leq s\leq t} \hatw^n(s)\right)\leq a_8 (t+1).
\ee
Consequently, $\barw^{n} = \frac{\hat W^n}{n^{1-\alpha}}\Go 0,$ and so from \eqref{diff-wl}
$\barq^n\Go 0$. Following the proof of Lemma \ref{conv} (i) (see equations \eqref{eqnA7}, \eqref{conv-A}, \eqref{eq:eq1333}, and \eqref{eqnA9}), we have $\bar{\T}^n\Go \bart^*$. \ink

\begin{lemma}\label{conv3}
Under the control policy $\{\T^n\}$, $(\hat{Q}^n_{\p(1)}, \ldots, \hat{Q}^n_{\p(K-1)})\Go 0$ as $n\go\infty.$
\end{lemma}
{\bf Proof:} Once again all the quantities in the proof  are considered under the control policy $\{\T^n\}$.
Without loss of generality, we assume $(\p(1), \ldots, \p(K))=(1,2,\ldots,K)$. For $t\geq 0$, let
$
\hatw_{-K}^n(t) := \sum_{j=1}^{K-1} \frac{\hatq^n_{j}(t)}{\mu^{n,*}_{j}}. 
$
It suffices to show that $\hatw^n_{-K}\Go 0$. Define for $t \ge 0$, $n \ge 1$,
 \bes\ba
 \hata^n_i(t) & = n^{-\alpha}\left(A^n_i(nt) - \int_0^{nt}\lambda^n_i(Y^n_s)ds\right) =  \hatn^n_{a,i}\left(\frac{1}{n}\int_0^{nt}\lambda^n_i(Y^n_s)ds\right), \\
  \hatd^n_i(t) & = n^{-\alpha}\left(D^n_i(nt)- \int_0^{nt}\mu^n_i(Y^n_s)\dot{\T}^n_{i}(s)ds\right) = \hatn^n_{s,i}\left(\frac{1}{n}\int_0^{nt}\mu^n_i(Y^n_s)\dot{\T}^n_i(s)ds\right).
 \ea\ees
From \eqref{queue-hat1} and \eqref{eqn-X}, we have for $t\geq 0$,
\bes\ba
 \hatw_{-K}^n(t)
& =\sum_{j=1}^{K-1} \frac{\hata^n_j(t)}{\mu^{n,*}_{j}}- \sum_{j=1}^{K-1}\frac{\hatd^n_j(t)}{\mu^{n,*}_{j}}+ n^{1-\alpha}\left(\sum_{j=1}^{K-1}\frac{\lambda_{j}^{n,*}}{\mu^{n,*}_{j}}-1\right)t + n^{-\alpha}\int_0^{nt} \left(1-\sum_{j=1}^{K-1} \dot{\T}^n_{j}(s)\right)ds\\
& + n^{-\alpha}\sum_{j=1}^{K-1} \frac{1}{ \mu^{n,*}_{j}}\int_0^{nt}\left(\lambda_{j}^n(Y^n_s)-\lambda^{n,*}_{j}\right)ds - n^{-\alpha} \sum_{j=1}^{K-1} \frac{1}{\mu^{n,*}_{j}}\int_0^{nt}\left(\mu_{j}^n(Y^n_s)-\mu^{n,*}_{j}\right)\dot{\T}^n_{j}(s)ds.
\ea\ees
Define, for $t\geq 0$,
\bes\ba
\xi^n(t) &= \sum_{j=1}^{K-1} \frac{\hata^n_j(t)}{\mu^{n,*}_{j}}- \sum_{j=1}^{K-1}\frac{\hatd^n_j(t)}{\mu^{n,*}_{j}} + n^{-\alpha} \sum_{j=1}^{K-1} \frac{1}{ \mu^{n,*}_{j}}\int_0^{nt}\left(\lambda_{j}^n(Y^n_s)-\lambda^{n,*}_{j}\right)ds \\
&\quad -  n^{-\alpha}\sum_{j=1}^{K-1} \frac{1}{\mu^{n,*}_{j}}\int_0^{nt}\left(\mu_{j}^n(Y^n_s)-\mu^{n,*}_{j}\right)\dot{\T}^n_{j}(s)ds,\\
\theta^n &=n^{1-\alpha}\left(\sum_{j=1}^{K-1}\frac{\lambda_{j}^{n,*}}{\mu^{n,*}_{j}}-1\right),\;\;
\zeta^n(t) = \frac{1}{n^\alpha}\int_0^{nt} \left(1-\sum_{j=1}^{K-1} \dot{\T}^n_{j}(s)\right)ds.
\ea\ees
Then for $t\geq 0$,
\bes
\hatw^n_{-K}(t) = \xi^n(t) + \theta^nt  + \zeta^n(t).
\ees
From Lemma \ref{conv}(ii), we have that $\xi^n$ converges to a Brownian motion with drift $0$ and variance $\sigma^2_{-K}$, where $\sigma^2_{-K}$ is given as follows:
\bes\ba
\mbox{\bf Case 1.}  \;\; & \sigma^2_{-K} = 2\sum_{j=1}^{K-1}\frac{\lambda^*_j}{(\mu^*_j)^2}. \\
\mbox{\bf Case 2.} \;\;  & \sigma^2_{-K} = 2\sum_{j=1}^{K-1}\frac{\lambda^*_j}{(\mu^*_j)^2} + \sum_{j=1}^{K-1}\frac{\Lambda_{jj}}{(\mu^*_j)^2}+ 2\sum_{1\leq i<j\leq K-1}\frac{\Lambda_{ij}}{\mu^*_i\mu^*_j}. \\
\mbox{\bf Case 3.} \;\; & \sigma^2_{-K} = \sum_{j=1}^{K-1}\frac{\Lambda_{jj}}{(\mu^*_j)^2}+ 2\sum_{1\leq i<j\leq K-1}\frac{\Lambda_{ij}}{\mu^*_i\mu^*_j}.
\ea\ees
Here $\Lambda$ is as defined in Lemma \ref{fclt-Y}. 
Next observe that
\bes
\zeta^n(t) = n^{1-\alpha}\int_0^t 1_{\{\hatq^n_{j}(s)=0, j=1,\ldots,K-1\}}ds = n^{1-\alpha}\int_0^t 1_{\{\hatw^n_{-K}(s)=0\}}ds,
\ees
which is nondecreasing and increases only when $\hatw^n_{-K}(\cdot)=0$. 
Define for $t\geq 0$,
\bes
\tau^n(t)= \sup\left\{0\leq s\leq t: \hatw^n_{-K}(s)=0\right\}.
\ees
Since $\hatw^n_{-K}(0)=0$, $\tau^n(t)\ge 0$, a.s. By definition,  $\hatw^n_{-K}(\tau^n(t)-)=0$. Also since $\hatw^n_{-K}(s)>0$ for all $s\in(\tau^n(t),t]$, we have $\zeta^n(\tau^n(t))=\zeta^n(t)$. Therefore, we have
\be\label{est7}\ba
0\leq \hatw^n_{-K}(t) &= \hatw^n_{-K}(t)-\hatw^n_{-K}(\tau^n(t)-) = \xi^n(t) - \xi^n(\tau^n(t)-)+ \theta^n(t-\tau^n(t)),
\ea\ee
and so
\be\label{est16}
-\frac{\theta^n}{n^{1-\alpha}}(t-\tau^n(t))\leq \frac{1}{n^{1-\alpha}}\left[\xi^n(t) - \xi^n(\tau^n(t)-)\right].
\ee
From the proof of Lemma \ref{conv4} (see \eqref{est-p1}, \eqref{est-p2}, \eqref{est-Phi}, and \eqref{est-G}), we have that 
\bes
\EE\left(\sup_{0\leq s\leq t} \frac{1}{n^{1-\alpha}}\left[|\xi^n(s)| + |\xi^n(\tau^n(s)-)|\right]\right)\go 0.
\ees
Also note that $$-\frac{\theta^n}{n^{1-\alpha}}\go -\left(\sum_{j=1}^{K-1}\frac{\lambda^*_{j}}{\mu^*_{j}}-1\right)>0.$$
Therefore, we have
\be\label{conv1}
\tau^n(t)\go t \;\; \mbox{in probability and uniformly on compact sets.}
\ee
Thus $\xi^n(t)-\xi^n(\tau^n(t)-)$ converges to $0$ in probability and uniformly on compact sets.
Using this observation in  \eqref{est16} we now have that 
$$\theta^n(\iota(\cdot) - \tau^n(\cdot))\Go 0,$$
where $\iota: [0,\infty)\go [0,\infty)$ is the identity map. 
Finally, from \eqref{est7}, we have
$\hatw^n_{-K}\Go 0.$
The result follows. \ink

{\bf Proof of Theorem \ref{optimal-limit}:} \\
As in the proof of Lemma \ref{conv3}, we assume
without loss of generality,  $(\p(1), \ldots, \p(K))=(1,2,\ldots,K)$. Also, once again  all the quantities in the proof are considered under control policy $\{\T^n\}$. From \eqref{est9} it follows that \eqref{eq:eq1556} holds with $\{\T^n\}$. Thus from Lemma \ref{conv}, we have $\hatx^n\Go X$, where $X$ is the Brownian motion introduced in Lemma \ref{conv}. 
Next recall that
$\hati^n$ as defined in \eqref{hati} satisfies: 
 $\hati^n(0)=0$; $\hati(\cdot)$ is nondecreasing; and increases only when $\hatw^n(\cdot)=0$. Therefore, from Definition \ref{sp-def}, we have for $t\geq 0$,
\be\label{workload1b}
\hatw^{n}(t)  = \Gamma\left(\sum_{i=1}^K\frac{\hatx^{n}_i}{\mu_i^{n,*}}\right)(t).
\ee
Using the Lipschitz continuity property of $\Gamma$ (Proposition \ref{sp}(ii)), and recalling $W^*$ defined in \eqref{opt-workload}, we have that 
\bes
\hatw^n \Go W^*.
\ees
Also  from Lemma \ref{conv3}, $(\hatq^n_1, \ldots, \hatq^n_{K-1})\Go 0$, and so from \eqref{diff-wl},
\be
\hatq^n_K\Go \mu^*_KW^*.
\ee 
To complete the proof we now establish a suitable uniform integrability estimate.
From \eqref{eqn-X} and \eqref{diff-wl}, we have 
\be\label{tilw}
\hatw^n(t) = \sum_{i=1}^K\frac{\hatq^n_i(t)}{\mu^{n,*}_i}=\sum_{i=1}^K\frac{\tilx^n_i(t)}{\mu^{n,*}_i} - n^{-\alpha + \frac{1}{2}}\sum_{i=1}^K\frac{\Phi^{\mu^n_i}(t)}{\mu^{n,*}_i}+ \hati^n(t),
\ee
where 
\bes\ba
\tilx^n_i(t) &= \hatn^n_{a,i}\left(\frac{1}{n}\int_0^{nt}\lambda^n_i(Y^n_s)ds\right)-\hatn^n_{s,i}\left(\frac{1}{n}\int_0^{nt}\mu^n_i(Y^n_s)\dott^n_i(s)ds\right) + n^{-(\alpha+\frac{\nu-1}{2})}G^{n^{\nu+1}}(t) \\
&\quad + \mu^{n,*}_in^{1-\alpha}\left(\frac{\lambda_i^{n,*}}{\mu^{n,*}_i}-\frac{\lambda_i^*}{\mu^*_i}\right)t.
\ea\ees
Next note that 
\be\label{est-I}\ba
\hati^n(t) &= -\inf_{0\leq s\leq t} \sum_{i=1}^K\frac{\hatx^n_i(t)}{\mu^{n,*}_i} = -\inf_{0\leq s\leq t} \left(\sum_{i=1}^K\frac{\tilx^n_i(t)}{\mu^{n,*}_i} - n^{-\alpha + \frac{1}{2}}\sum_{i=1}^K\frac{\Phi^{\mu^n_i}(t)}{\mu^{n,*}_i} \right)\\
& \leq -\inf_{0\leq s\leq t} \sum_{i=1}^K\frac{\tilx^n_i(t)}{\mu^{n,*}_i} - n^{-\alpha + \frac{1}{2}}\inf_{0\leq s\leq t} \sum_{i=1}^K\frac{- \Phi^{\mu^n_i}(t)}{\mu^{n,*}_i},
\ea\ee
and
\be\label{est-I1}\ba
\hati^n(t) & \geq -\sup_{0\leq s\leq t} \sum_{i=1}^K\frac{|\tilx^n_i(t)|}{\mu^{n,*}_i} - n^{-\alpha + \frac{1}{2}}\inf_{0\leq s\leq t} \sum_{i=1}^K\frac{- \Phi^{\mu^n_i}(t)}{\mu^{n,*}_i},
\ea\ee
Define for $t\geq 0$,
\bes
\tilq^n_j(t) = \hatq^n_j(t), \;\; j =1, 2, \ldots, K-1,
\ees
and 
\bes
\tilq^n_K(t) = \hatq^n_K(t) +  \mu^{n,*}_K n^{-\alpha + \frac{1}{2}}\sum_{i=1}^K\frac{\Phi^{\mu^n_i}(t)}{\mu^{n,*}_i} + \mu^{n,*}_K n^{-\alpha + \frac{1}{2}}\inf_{0\leq s\leq t}\sum_{i=1}^K\frac{-\Phi^{\mu^n_i}(s)}{\mu^{n,*}_i}.
\ees
Let for $t\geq 0$,
\bes\ba
\tilw^n(t) &= \sum_{i=1}^K\frac{\tilq^n_i(t)}{\mu^{n,*}_i} = \sum_{i=1}^K\frac{\hatq^n_i(t)}{\mu^{n,*}_i} + n^{-\alpha + \frac{1}{2}}\sum_{i=1}^K\frac{\Phi^{\mu^n_i}(t)}{\mu^{n,*}_i} + n^{-\alpha + \frac{1}{2}} \inf_{0\leq s\leq t}\sum_{i=1}^K\frac{-\Phi^{\mu^n_i}(s)}{\mu^{n,*}_i} \\
& = \sum_{i=1}^K\frac{\tilx^n_i(t)}{\mu^{n,*}_i} + \hati^n(t) + n^{-\alpha + \frac{1}{2}} \inf_{0\leq s\leq t}\sum_{i=1}^K\frac{-\Phi^{\mu^n_i}(s)}{\mu^{n,*}_i}.
\ea\ees
where the last equality uses \eqref{tilw}.
Using the estimates in \eqref{est-I} and \eqref{est-I1}, we observe that
\be\label{eq:eq1600}\ba
 &\sum_{i=1}^K\frac{\tilx^n_i(t)}{\mu^{n,*}_i}-\sup_{0\leq s\leq t} \sum_{i=1}^K\frac{|\tilx^n_i(t)|}{\mu^{n,*}_i}\leq \tilw^n(t) \leq \sum_{i=1}^K\frac{\tilx^n_i(t)}{\mu^{n,*}_i} -\inf_{0\leq s\leq t} \sum_{i=1}^K\frac{\tilx^n_i(t)}{\mu^{n,*}_i}= \Gamma\left(\sum_{i=1}^K\frac{\tilx^n_i(t)}{\mu^{n,*}_i}\right).
\ea\ee
 From \eqref{est-p1} -- \eqref{est-drift}, and \eqref{est-G}, we have 
for some $a_1 \in (0, \infty)$, 
$$
\sup_n \EE\left(\sup_{0\leq s\leq t} |\tilx^n(s)|^2\right) \leq a_1 (t+1), \mbox{ for all } t \ge 0.
$$
Using this estimate in \eqref{eq:eq1600} we have,
for some $a_2\in (0,\infty)$, 
\be\label{UI}
\EE\left(\sup_{0\leq s\leq t} \tilw^n(s)\right)^2 \leq a_2 (t+1). 
\ee
From Lemma \ref{ergo},   $n^{-\alpha + \frac{1}{2}}\Phi^{\mu^n_i}\Go 0$, and so
\bes
\tilq^n_K \Go \mu^*_KW^*, \;\; \tilq^n_j\Go 0, \; j=1,2,\ldots, K-1.
\ees
By \eqref{UI}, we now have that
\bes
\EE\left(\int_0^\infty e^{-\gamma t}( c\cdot \tilq^n(t)) dt \right)\go \EE\left(\int_0^\infty e^{-\gamma t} c_K\mu^*_KW^*(t) dt \right) = J^*.
\ees
Since $\tilde Q^n = \hat Q^n$ in cases 1(b), 2 and 3, we have the result for these cases.  Finally consider the case 1(a).
Using Lemma \ref{ergo}, we have for some $a_3, a_4 \in(0,\infty)$,
\bes\ba
&\left|\EE\left(\int_0^\infty e^{-\gamma t} c\cdot \hatq^n(t) dt \right) - \EE\left(\int_0^\infty e^{-\gamma t} c\cdot \tilq^n(t) dt \right)\right| \leq \EE\left(\int_0^\infty e^{-\gamma t} c\cdot |\hatq^n(t)-\tilq^n(t)| dt \right)\\
&\leq a_3\sum_{i=1}^K\int_0^\infty e^{-\gamma t} \EE\left(\sup_{0\leq s\leq t} |\Phi^{\mu^n_i}(s)| \right) dt \leq  \frac{a_4}{n^{\nu-1/2}} \int_0^\infty e^{-\gamma t}  t dt \go 0. 
\ea\ees
The result follows.  \ink

{\bf Proof of Theorem \ref{optimal}:} Let $\{T^n\}$ be a sequence of admissible control policies. 
The result holds trivially if $\underline{J}(\{T^n\}) = \liminf_{n\to \infty} \hat J^n(T^n) = \infty$.
Suppose now that
$\underline{J}(\{T^n\})<\infty$. 
Let $\{n'\}$ be a subsequence such that $\hat J^{n'}(T^{n'}) \to \underline{J}\{T^n\}$, as $n' \to \infty$.  Henceforth we relabel the subsequence
$\{n'\}$ as $\{n\}$.
Throughout the proof, all quantities are considered under the control sequence
$\{T^n\}$. Recalling \eqref{diff-wl} and that
$\hat W^n(t) \ge 0$; $\hati^n$ is nondecreasing and starts from zero, we have from Proposition \ref{sp} (iii)  that for $t\geq 0$,
\be\label{min}
\hati^n(t) \geq \hati^{n.*}(t):=-\inf_{0\leq s\leq t} \sum_{i=1}^K\frac{\hatx^{n}_i(s)}{\mu_i^{n,*}},
\ee
and with 
$\hatw^{n,*}(t):= \sum_{i=1}^K\frac{\hatx^{n}_i(t)}{\mu_i^{n,*}} + \hati^{n,*}(t)$
we have for all $t\geq 0$,
\be\label{min1}
\hatw^n(t)\geq \hatw^{n,*}(t).
\ee
Also from Definition \ref{sp-def}, we have 
$
\hatw^{n,*}(t)=\Gamma\left(\sum_{i=1}^K\frac{\hatx^{n}_i}{\mu_i^{n,*}}\right)(t).
$ 
From Lemma \ref{conv}(ii), and using the Lipschitz continuity of $\Gamma$ (Proposition \ref{sp} (ii)), we now have 
\be\label{workload-conv}
\hatw^{n,*}\Go W^*,
\ee
where $W^*$ is defined by \eqref{opt-workload}. 

Next, for fixed $w\geq 0$  consider the following linear programming (LP) problem:
\be\label{LP}\ba
\mbox{Minimize}_{q\in\RR^K_+} \  \; \sum_{i=1}^Kc_iq_i \;\;
\quad \mbox{subject to} \; \sum_{i=1}^K\frac{q_i}{\mu^*_i} = w.
\ea\ee
Elementary calculations show that the the value of the above LP problem is
\be \label{eq:eq1657}
V(w)= c_{\p(K)}\mu^*_{\p(K)} w,
\ee
and the corresponding solution is
\be
q_{\p(K)} = \mu_{\p(K)}^* w, \; \mbox{and} \; q_j =0, j= \p(1),\ldots,\p(K-1).
\ee
Therefore, letting $\tilde W^n(t) = \sum_{i=1}^K \frac{\hat Q^n_i(t)}{\mu_i^*}$,
for $c \ge 0$,
$$
c\cdot \hat Q^n(t) \ge V(\tilde W^n(t)) = V(\hat W^n(t)) + c_K\mu^*_K \sum_{i=1}^K \hat Q^n_i(t)\left[\frac{1}{ \mu^{n,*}_i} - \frac{1}{ \mu^{*}_i}\right]
\ge V(\hat W^n(t)) + Z^n(t).$$
where once again, without loss of generality we assume,  $(\p(1), \ldots, \p(K))=(1,2,\ldots,K)$ and
$Z^n(t) = c_K\mu^*_K \sum_{i=1}^K \hat Q^n_i(t)[\frac{1}{ \mu^{n,*}_i} - \frac{1}{ \mu^{*}_i}]$.
Since $\underline{J}(\{T^n\})<\infty$,
$$
\limsup_{n\to \infty} \EE \int_0^{\infty} e^{-\gamma t} |Z^n(t)| dt \le 
\frac{c_K\mu^*_K}{\underline{c}}
\lim_{n\to \infty} \max_i\left|\frac{1}{ \mu^{n,*}_i} - \frac{1}{ \mu^{*}_i}\right|
\hat J^n(T^n) = 0,$$
where $\underline{c} = \min_i\{c_i\}$.
Thus from \eqref{workload-conv}
\bes\ba
\lim_{n\to\infty} \hat J^n(T^n) = \underline{J}(\{T^n\})&\geq 
\liminf_{n\to \infty}\EE\left(\int_0^\infty e^{-\gamma t} V(\hat W^{n,*}(t)) dt\right)\\
&\geq\EE\left(\int_0^\infty e^{-\gamma t}V(W^{*}(t))dt\right)\\
&= J^* = \lim_{n\go\infty}\hatj^n(\T^n),
\ea\ees
where the last equality follows from Theorem \ref{optimal-limit}.
\ink


\setcounter{equation}{0}
\appendix
\numberwithin{equation}{section}

\section{Proof of Lemma \ref{bcp-soln}}\label{bcp2}
We begin by introducing a $1$-dimensional control problem referred to as the `workload control problem'
\begin{definition} The workload control problem (WCP) is to find an $\RR_+$-valued RCLL $\{\bar \clf_t\}_{t\ge 0}$ adapted stochastic process $\tili$ which minimizes
\bes
\EE\left(\int_0^\infty e^{-\gamma t} V(\tilw(t))dt\right),
\ees
subject to
\be\label{ewf-cond}\ba
&\tilw(t) = \sum_{i=1}^K\frac{X_i(t)}{\mu^*_i}+ \tili(t)\geq 0, \; \mbox{for all} \; t\geq 0, \\
&\tili(\cdot) \; \mbox{is nondecreasing and} \; \tili(0)=0,
\ea\ee
where $X$ and $\{\bar \clf_t\}_{t\ge 0}$ are as in Definition \ref{bcp-def} and $V$ is defined in \eqref{eq:eq1657}.
\end{definition}

{\bf Proof of Lemma \ref{bcp-soln}:} 
Recall $(W^*, I^*)$, defined in \eqref{opt-workload}.
From  \eqref{optimal-q} we see that, for every $t\ge 0$, $Q^*(t)$ is a solution of the LP problem \eqref{LP} associated with $W^*(t)$. Thus 
\be\label{eq:eq1651}
\sum_{i=1}^Kc_iQ^*_i(t)=V(W^*(t)).
\ee
Also note that the pair $(W^*, I^*)$ satisfies the two conditions in \eqref{ewf-cond}.  Furthermore, if
$(\tilw, \tili)$ is another pair satisfying \eqref{ewf-cond}, we have from Proposition \ref{sp}(iii) that
$\tilw(t) \ge W^*(t)$ and consequently, 
$
V(\tilw(t)) \ge V(W^*(t))$,  for all $t \ge 0$, where $V$ is as defined in \eqref{eq:eq1657}.

Now let $\tilde \eta \in \tilde \cla$ and define the corresponding $\tilde Q_i$, $i = 1, \cdots K$, and $\tilde I$ through \eqref{eqn8}.
Define 
\be\label{eq:eq1702}
\tilw(t) = \sum_{i=1}^K \frac{\tilq_i(t)}{\mu^*_i}.
\ee  Clearly $(\tilw, \tili)$  satisfy \eqref{ewf-cond} and so from the
above discussion 
\be
\label{eq:eq1756}
V(\tilde W(t)) \ge V(W^*(t)), \mbox{ for all } t \ge 0. \ee
Also, from the definition of $V$ and \eqref{eq:eq1702} we see that
$\sum_{i=1}^K c_i \tilde Q_i(t) \ge V(\tilw(t))$, for all $t \ge 0$.  Combining this inequality with \eqref{eq:eq1651} and \eqref{eq:eq1756}
we get that $\sum_{i=1}^K c_i \tilde Q_i(t) \ge \sum_{i=1}^K c_i \tilde Q_i^*(t)$, for all $t \ge 0$.  Since $\tilde \eta \in \tilde \cla$
is arbitrary, we have the result.
\ink

\section{Proof of \eqref{fclt-m}} \label{lemma}
We only treat the case when $K=1$.  The general case can be proved similarly.
We first consider the special case  when $\tily^n$ is stationary. Denote by $\EE_{\pi^n}$ the expectation operator when
$\tily^n(0)$ has distribution $\pi^n$.
Let $\m^n(t) = \parvs(n)^{-\frac{1}{2}}M^n_{\parvs(n)t}, t\geq 0$. Noting that $|\Delta \m^n_t|\leq a_1\parvs(n)^{-\frac{1}{2}}$ for all $t\geq 0$, we have that for any $\epsilon>0$ and $t>0$,
$$\lim_{n\to\infty}\EE_{\pi^n}\left(\sum_{0\leq s\leq t} |\Delta \m^n_t|^21_{\{|\Delta \m^n_t|>\epsilon\}}\right) = 0.$$
By martingale central limit theorem (see Corollary VIII.3.24 in \cite{jacod}), it suffices to show that
$ \langle \m^n\rangle_t \to - 2 t \pi(\hat\lambda\clq\hat\lambda)$ in probability for all $t\geq 0$, where $ \langle \m^n\rangle$ denotes the predictable quadratic variation
of the martingale $\m^n$. We note from \eqref{martingale} that
\bes
\m^n(t) = \frac{1}{\sqrt{\betan}} \left(\hat\lambda^n(\tily^n({\betan t})) - \hat\lambda^n(\tily^n(0)) -\int_0^{\betan t} \clq^n \hat\lambda^n(\tily^n_s) ds\right).
\ees
Using Lemma VIII.3.68 in \cite{jacod}, we have 
\bes
 \langle \m^n\rangle_t = \frac{1}{\betan}\int_0^{\betan t} \left( \clq^n(\hat\lambda^n)^2(\tily^n_s) - 2\hat\lambda^n \clq^n\hat\lambda^n(\tily^n_s) \right)ds.
\ees
We then note that
\bes\ba
& \EE_{\pi^n}\left( \langle \m^n\rangle_t  + 2 \pi(\hat\lambda\clq\hat\lambda)\right)^2 \\
& =  \EE_{\pi^n}\left[ \frac{1}{\betan} \sum_{1\leq k\leq \nt} \int_{k-1}^k\left(\clq^n(\hat\lambda^n)^2(\tily^n_s) - 2\hat\lambda^n \clq^n\hat\lambda^n(\tily^n_s)  + 2 \pi(\hat\lambda\clq\hat\lambda) \right)ds \right]^2 + \frac{a_2}{\betan^2}\\
& \leq \frac{a_2}{\betan^2}+\frac{1}{\betan^2} \sum_{1\leq k\leq \nt} \EE_{\pi^n}\left[ \int_{k-1}^k\left(\clq^n(\hat\lambda^n)^2(\tily^n_s) - 2\hat\lambda^n \clq^n\hat\lambda^n(\tily^n_s)  + 2 \pi(\hat\lambda\clq\hat\lambda) \right)ds\right]^2 \\
&\quad + \frac{2}{\betan^2} \sum_{1\leq k_1<k_2\leq \nt} \EE_{\pi^n}\left[\int_{k_1-1}^{k_1}\left(\clq^n(\hat\lambda^n)^2(\tily^n_s) - 2\hat\lambda^n \clq^n\hat\lambda^n(\tily^n_s)  + 2 \pi(\hat\lambda\clq\hat\lambda) \right)ds \right. \\
&\qquad \times \left. \int_{k_2-1}^{k_2}\left(\clq^n(\hat\lambda^n)^2(\tily^n_s) - 2\hat\lambda^n \clq^n\hat\lambda^n(\tily^n_s)  + 2 \pi(\hat\lambda\clq\hat\lambda) \right)ds\right] \\
& \leq \frac{a_3}{\betan}+ \frac{2}{\betan^2}  \sum_{1\leq k_1<k_2\leq \nt}\EE_{\pi^n}\left[\int_{k_1-1}^{k_1}\left(\clq^n(\hat\lambda^n)^2(\tily^n_s) - 2\hat\lambda^n \clq^n\hat\lambda^n(\tily^n_s)  + 2 \pi(\hat\lambda\clq\hat\lambda) \right)ds \right. \\
&\qquad \times \left. \EE_{\pi^n}\left(\left.\int_{k_2-1}^{k_2}\left(\clq^n(\hat\lambda^n)^2(\tily^n_s) - 2\hat\lambda^n \clq^n\hat\lambda^n(\tily^n_s)  + 2 \pi(\hat\lambda\clq\hat\lambda) \right)ds\right| \clg^n_{k_1}\right)\right] \\
& \leq  \frac{a_3}{\betan} + \frac{2}{\betan^2}  \sum_{1\leq k_1<k_2\leq \nt}\EE_{\pi^n}\left[\int_{k_1-1}^{k_1}\left(\clq^n(\hat\lambda^n)^2(\tily^n_s) - 2\hat\lambda^n \clq^n\hat\lambda^n(\tily^n_s)  + 2 \pi(\hat\lambda\clq\hat\lambda) \right)ds \right. \\
&\qquad \times \left. \int_{k_2-1}^{k_2}\left(P^{n,s-k_1}(\tily^n_{k_2-1}, \clq^n(\hat\lambda^n)^2) - 2P^{n,s-k_1}(\tily^n_{k_2-1}, \hat\lambda^n \clq^n\hat\lambda^n)  + 2 \pi(\hat\lambda\clq\hat\lambda) 
\right)ds \right]
\ea\ees
Noting that $\pi^n(\clq^n(\hat\lambda^n)^2) = 0$ and $\pi^n(\hat\lambda^n\clq^n\hat\lambda^n)\to\pi(\hat\lambda\clq\hat\lambda)$, and using Assumption \ref{uniferg}, we have that for $k_2-1\leq s \leq k_2,$
\bes\ba
\sup_{y\in\LL}\left|P^{n,s-k_1}(y, \clq^n(\hat\lambda^n)^2)- 2P^{n,s-k_1}(y, \hat\lambda^n \clq^n\hat\lambda^n)  + 2 \pi(\hat\lambda\clq\hat\lambda)\right| \leq a_4 \rho^{s-k_1}.
\ea\ees
Thus we have that 
\bes\ba
 \EE\left( \langle \m^n\rangle_t  + 2 \pi(\hat\lambda\clq\hat\lambda)\right)^2 & \leq  \frac{a_3}{\betan} + \frac{a_5}{\betan^2}  \sum_{1\leq k_1<k_2\leq \nt}\int_{k_2-1}^{k_2} \rho^{s-k_1}ds\\
 & \to 0 \ \mbox{as} \ n\to\infty.  
\ea\ees
This completes the proof for the case when $\tily^n$ is stationary. 

Finally we consider the case when $\tily^n(0)$ has an arbitrary  distribution. It suffices to consider the setting where $\tily^n(0) = x \in \LL$,
for all $n \ge 1$. The argument below follows the proof of Theorem 4.3 in \cite{glynn-meyn96}. 
Define for $r\in (0,\infty)$,
\bes
G^{n,r}(t) = \frac{1}{\sqrt{\betan}} \int_r^{\betan t+r} (\lambda^n(\tily^n_s)-\lambda^{n,*})ds.
\ees
Clearly, for any $r\in (0,\infty)$,
\be\label{shift}
\sup_{t\geq 0}|G^n(t) - G^{n,r}(t)| \go 0 \;\; \mbox{as} \; n\go \infty. 
\ee
Let $\phi$ be a real valued bounded Lipschitz function on $\mathcal{C}([0, \infty), \RR)$. Then from \eqref{shift}, we have, as $n \to \infty$,
\bes
|\EE_x[\phi(G^n)] - \EE_x[\phi(G^{n,r})]| \go 0,
\ees
where $\EE_x$ denotes the expectation operator when $\tily^n(0) =x$ a.s.
Using the Markov property, we have
\be \label{eq:eq1822}
\EE_x[\phi(G^{n,r})] = \int_{\LL} P^{n,r}(x, dy) \EE_y[\phi(G^n)]. 
\ee
Thus we have for some $a_6\in (0,\infty)$,
\bes\ba
& \left|\EE_x[\phi(G^n)]-  \EE_{\pi^n}[\phi(G^n)] \right|\\
& = \left|\EE_x[\phi(G^n)] - \EE_x[\phi(G^{n,r})] + \int_{\LL} (P^{n,r}(x, dy)-\pi^n(y)) \EE_y[\phi(G^n)]\right| \\
&\leq |\EE_x[\phi(G^n)] - \EE_x[\phi(G^{n,r})]| + a_6 \|P^{n,r}(x, \cdot)- \pi^n(\cdot)\|
\ea\ees
Using \eqref{eq:eq1822} and Assumption \eqref{uniferg} we see that the last expression converges to $0$ on sending first
 $n\go \infty$ and then  $r\go\infty$.
This completes the proof of \eqref{fclt-m}.   \ink

\section{Simulation results in Example \ref{ex}} \label{simulation}

We continue with Example \ref{ex}. We plot the total cost and total discounted cost at each discrete time. More precisely, recall that the holding cost $c=(20,25)'$ and discount factor $\gamma=2$. For $10$ sample paths, define
\bes
C_1^n(t) = \frac{1}{10}  \sum_{s=1}^{10}\left(c_1q^n_{1,s}(t) + c_2q^n_{2,s}(t)\right),
\ees
and
\bes
C_2^n(t) = \frac{1}{10}  \sum_{s=1}^{10}\left(c_1e^{-\gamma t} q^n_{1,s}(t) + c_2e^{\gamma t}q^n_{2,s}(t)\right). 
\ees
Here $C_1^n(t)$ measures the total cost for the queueing system at time $t$, while $C_2^n(t)$ gives the discount cost at time $t$. 

Let $n=100$. We plot $C^n_1(t), t=0.1, 0.5, 1.0, \ldots, 100$ and $C^n_2(t), t=0.1, 0.2, 0.3, \ldots, 10$ for six different scalings. In all the following figures, the first plot gives $C^n_1$ and the second one is $C^n_2$. We use red color to represent $c\mu^*$-rule and green color for dynamic $c\mu$-rule. We observe that in all the plots, the red line is located below the green one for most of time, which says that the cost under $c\mu^*$-rule is smaller. 

\vspace{1in}

\begin{figure}[htbp]
\includegraphics[width=160mm, height=45mm]{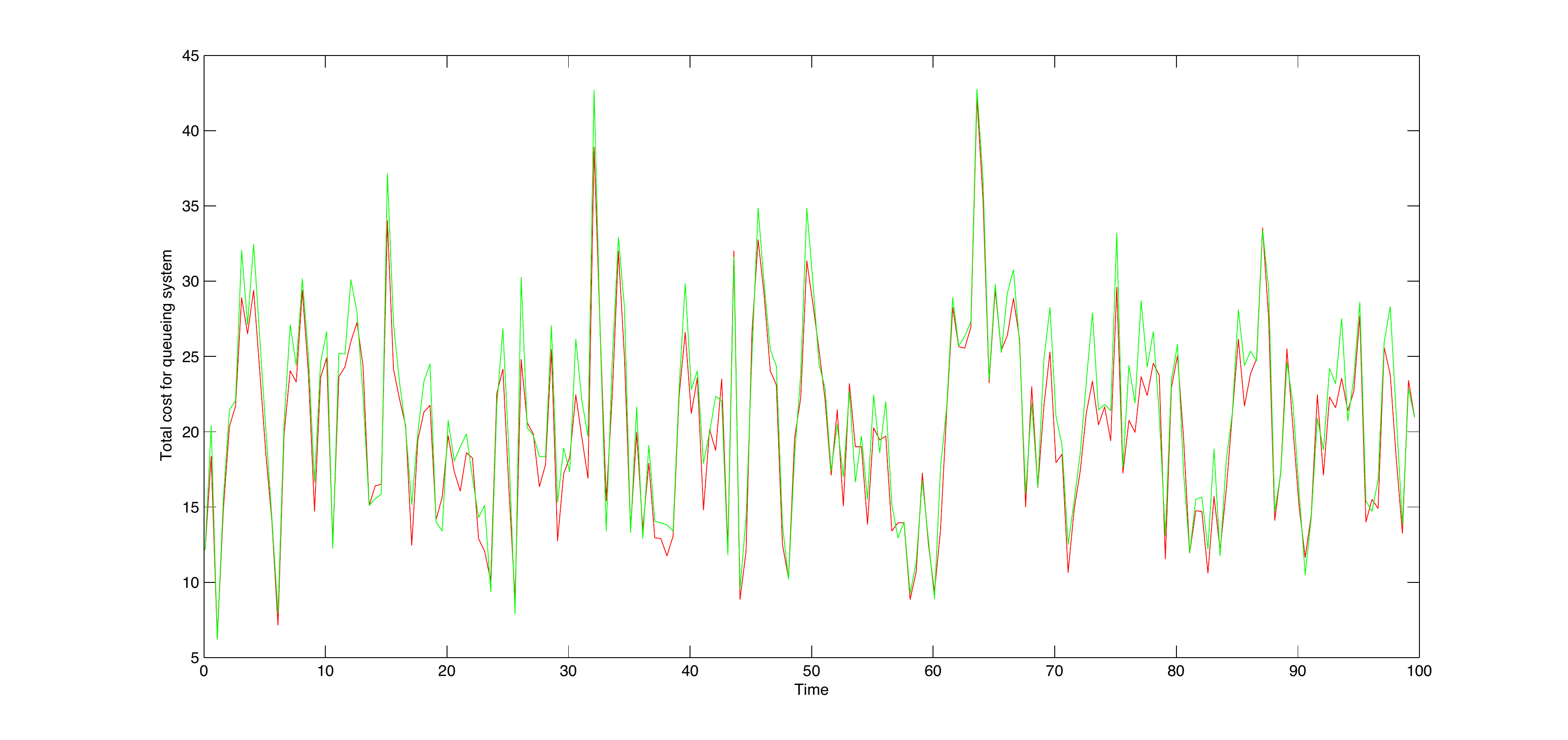}
\includegraphics[width=160mm, height=45mm]{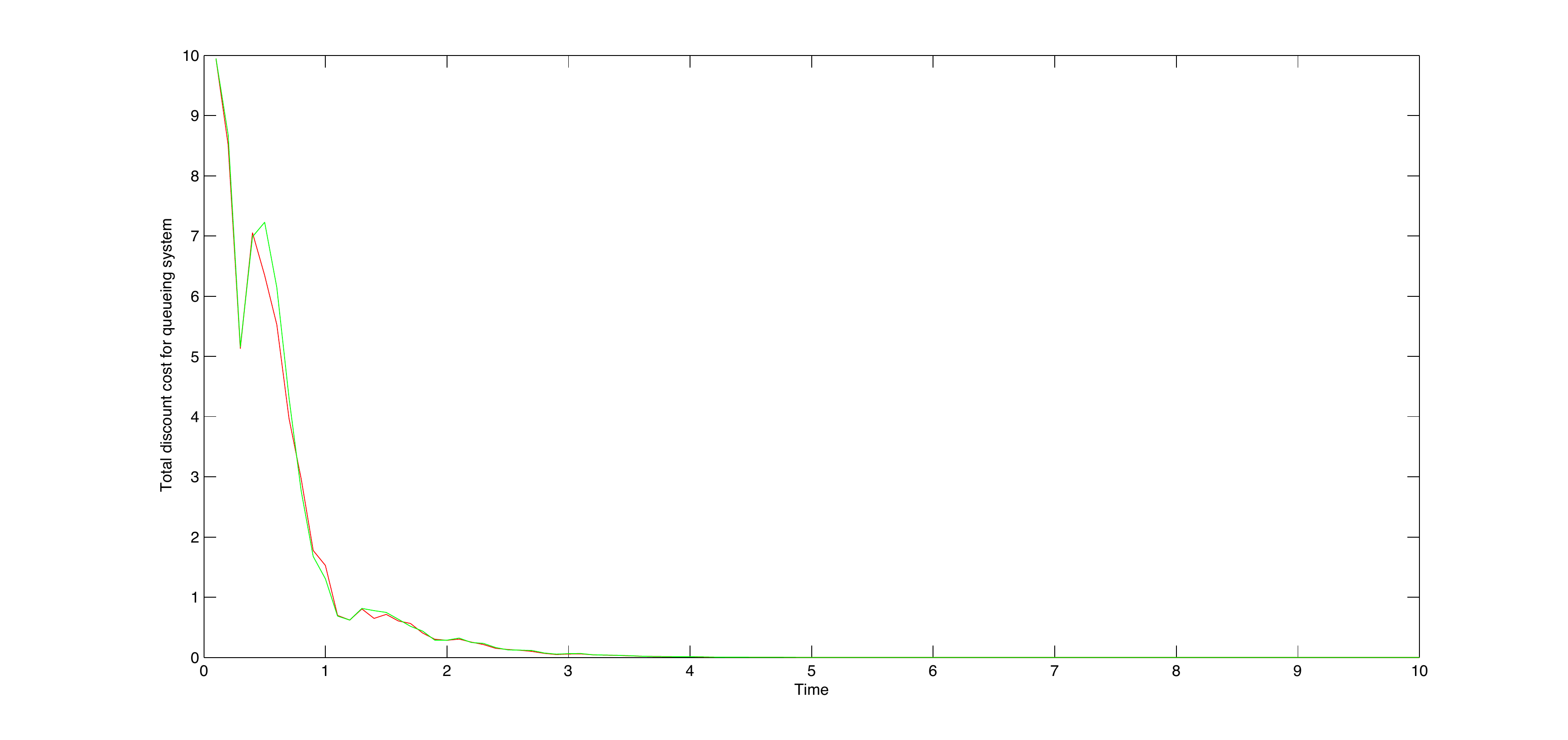}
\caption{Average sample path for the case when $\nu= 1$ and $\alpha=1/2$. }
\end{figure}

\begin{figure}[htbp]
\includegraphics[width=160mm, height=45mm]{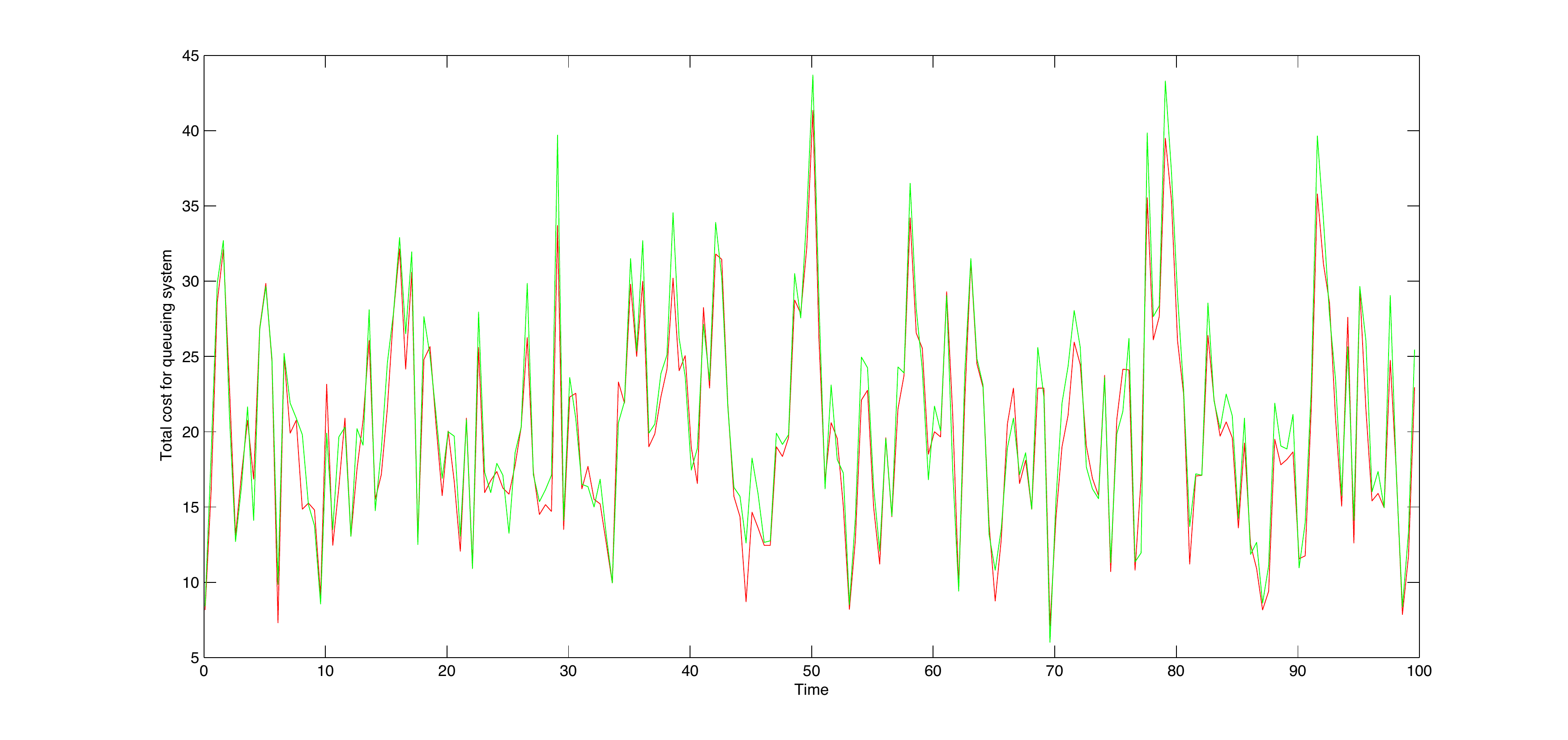}
\includegraphics[width=160mm, height=45mm]{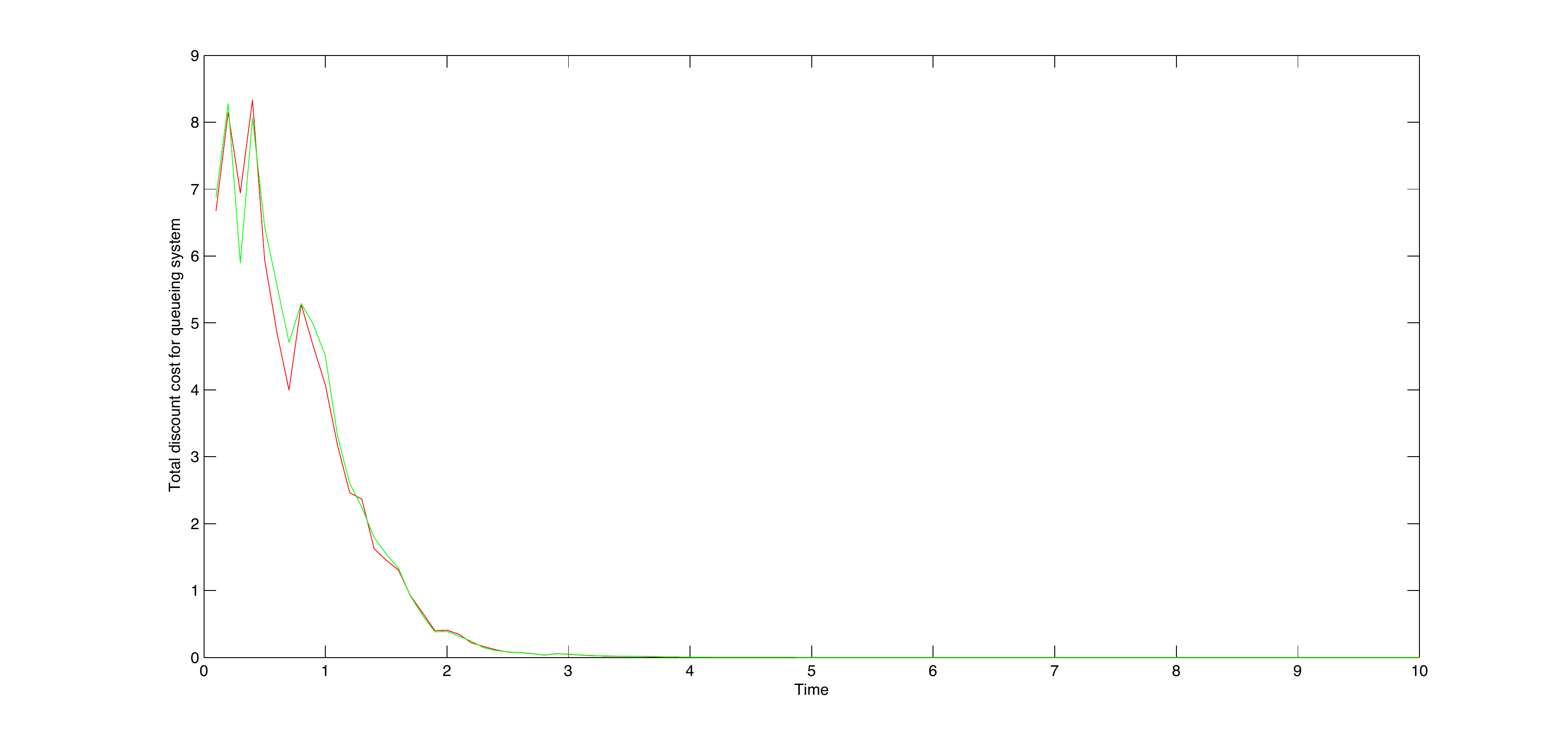}
\caption{Average sample path for the case when $\nu= 2/3$ and $\alpha=1/2$. }
\end{figure}

\begin{figure}[htbp]
\includegraphics[width=160mm, height=45mm]{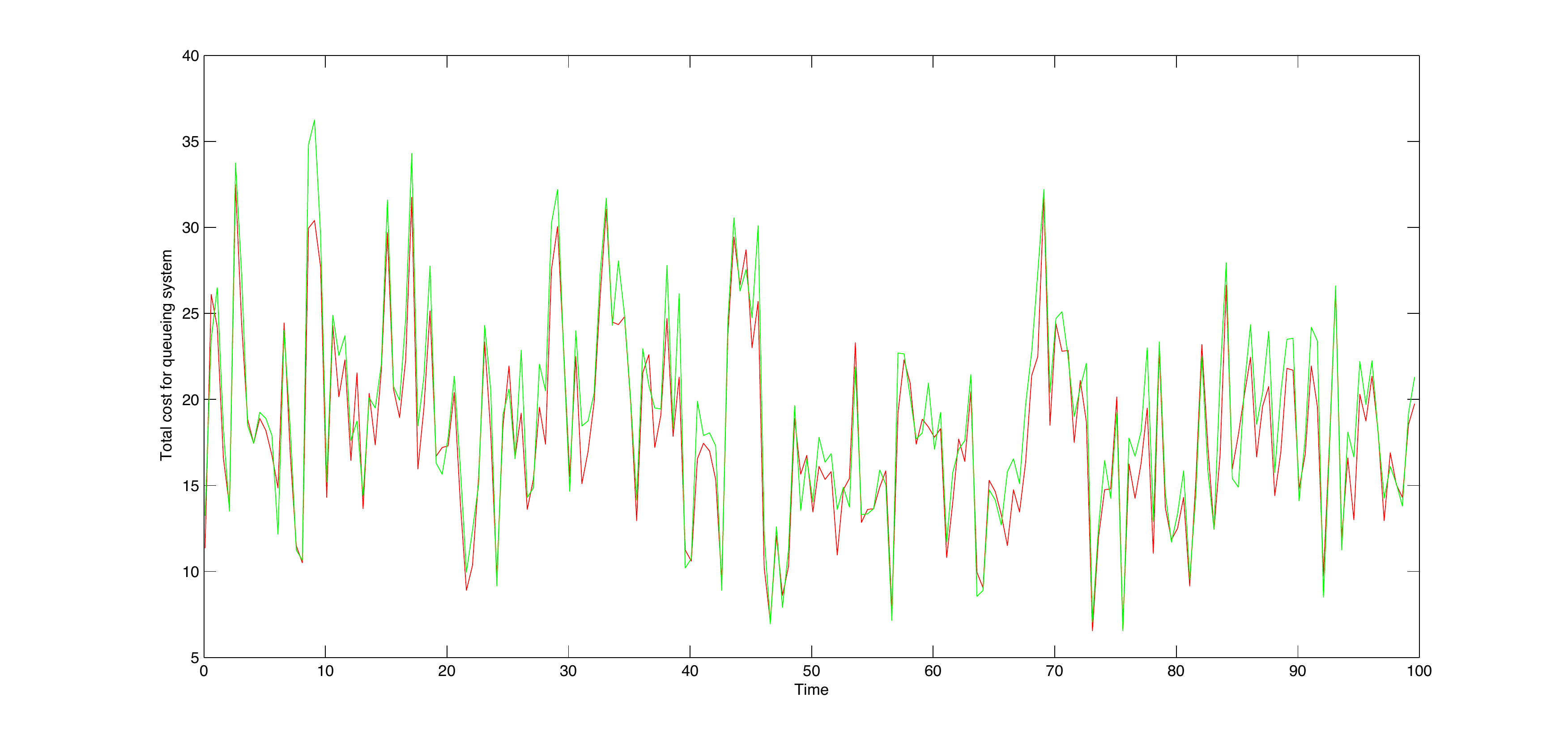}
\includegraphics[width=160mm, height=45mm]{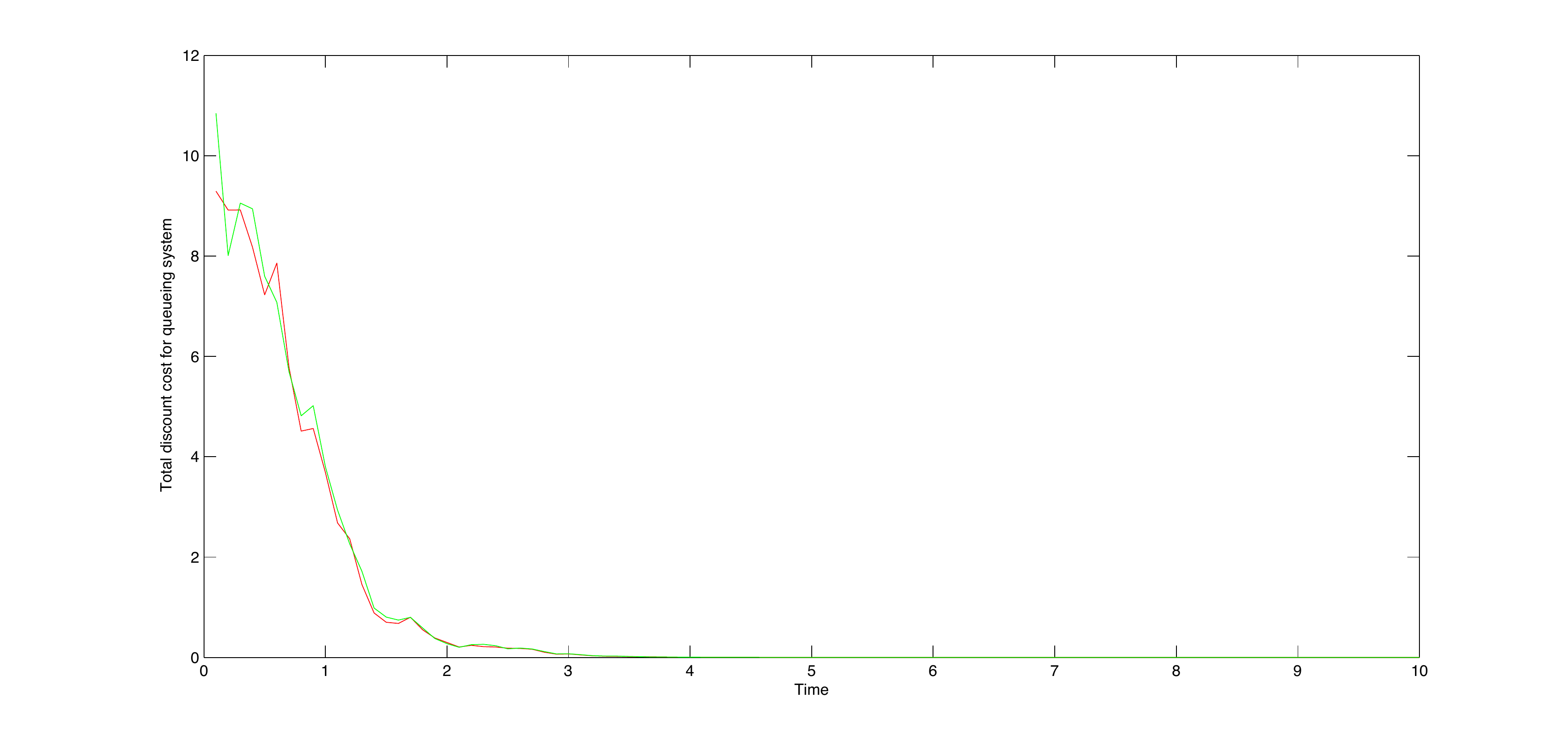}
\caption{Average sample path for the case when $\nu= 1/3$ and $\alpha=1/2$. }
\end{figure}

\begin{figure}[htbp]
\includegraphics[width=160mm, height=45mm]{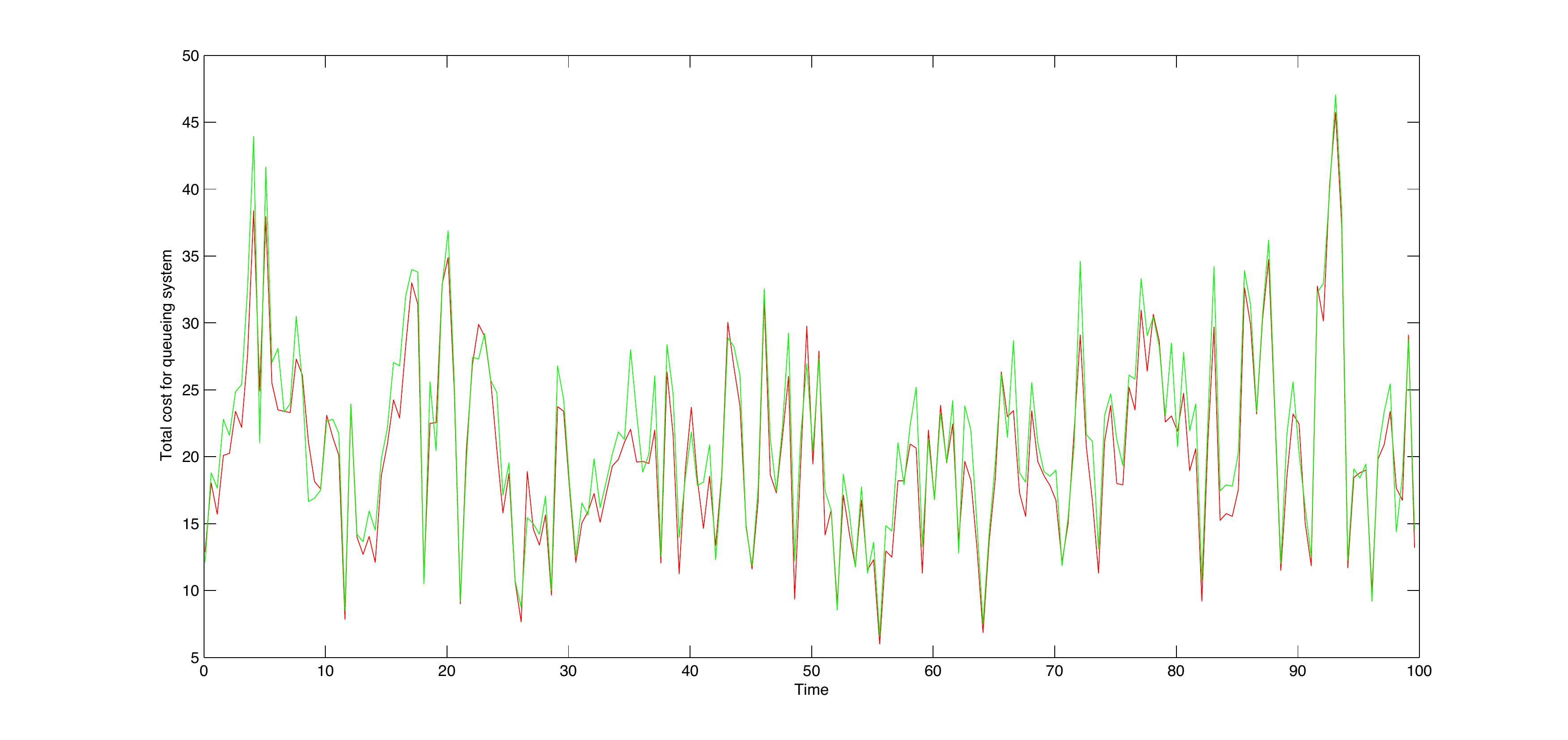}
\includegraphics[width=160mm, height=45mm]{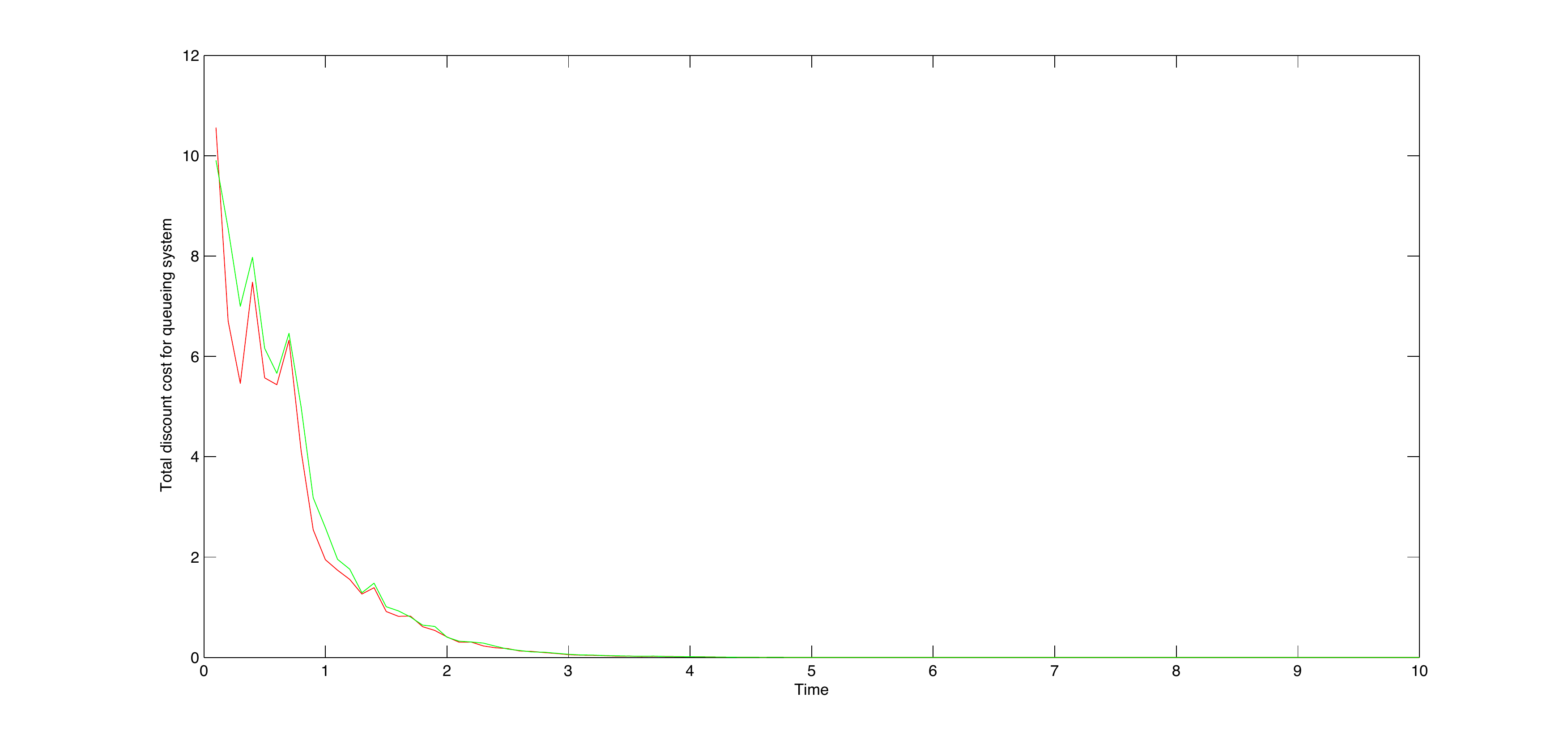}
\caption{Average sample path for the case when $\nu= 0$ and $\alpha=1/2$ }
\end{figure}

\begin{figure}[htbp]
\includegraphics[width=160mm, height=45mm]{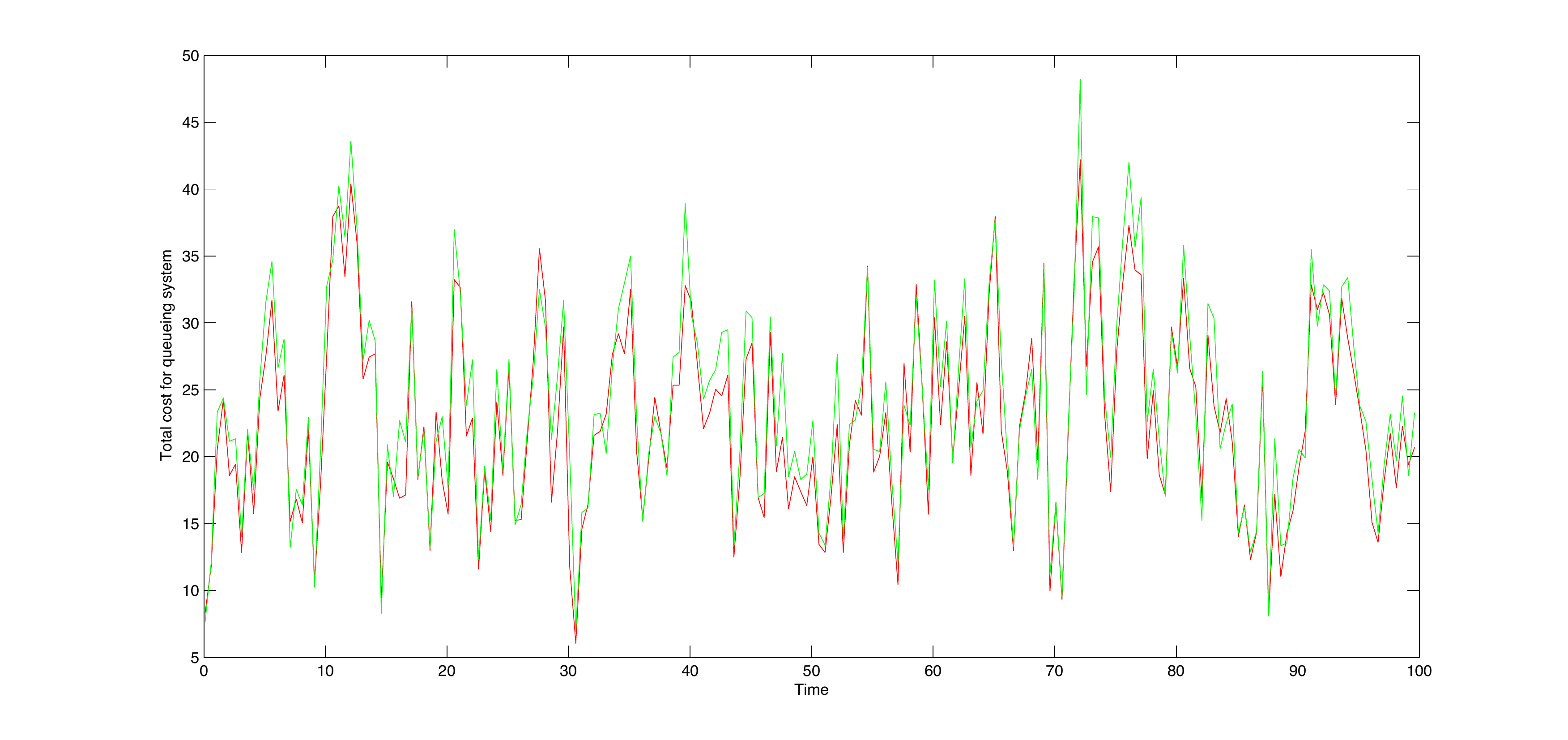}
\includegraphics[width=160mm, height=45mm]{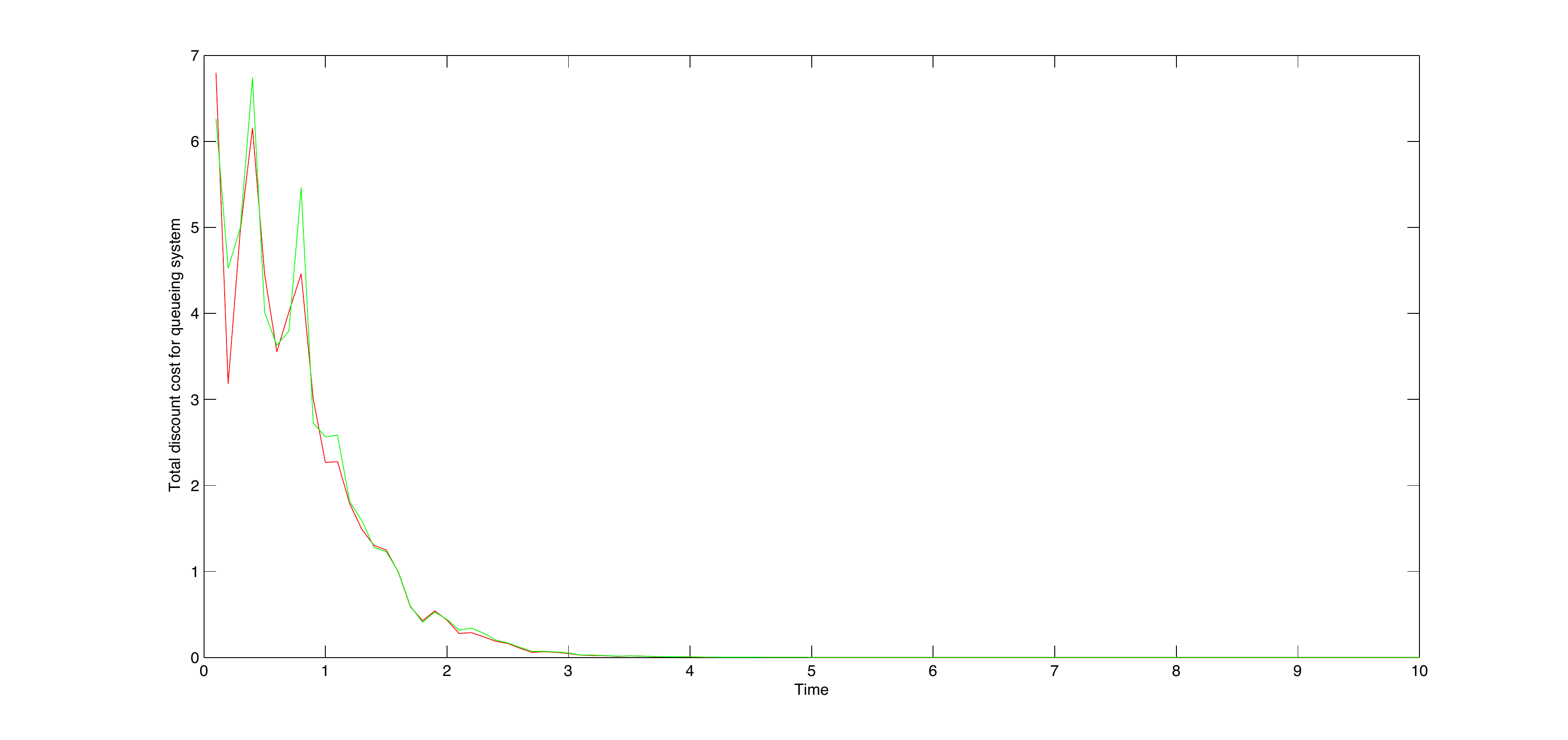}
\caption{Average sample path for the case when $\nu= -1/3$ and $\alpha=1/2$. }
\end{figure}

\begin{figure}[htbp]
\includegraphics[width=160mm, height=45mm]{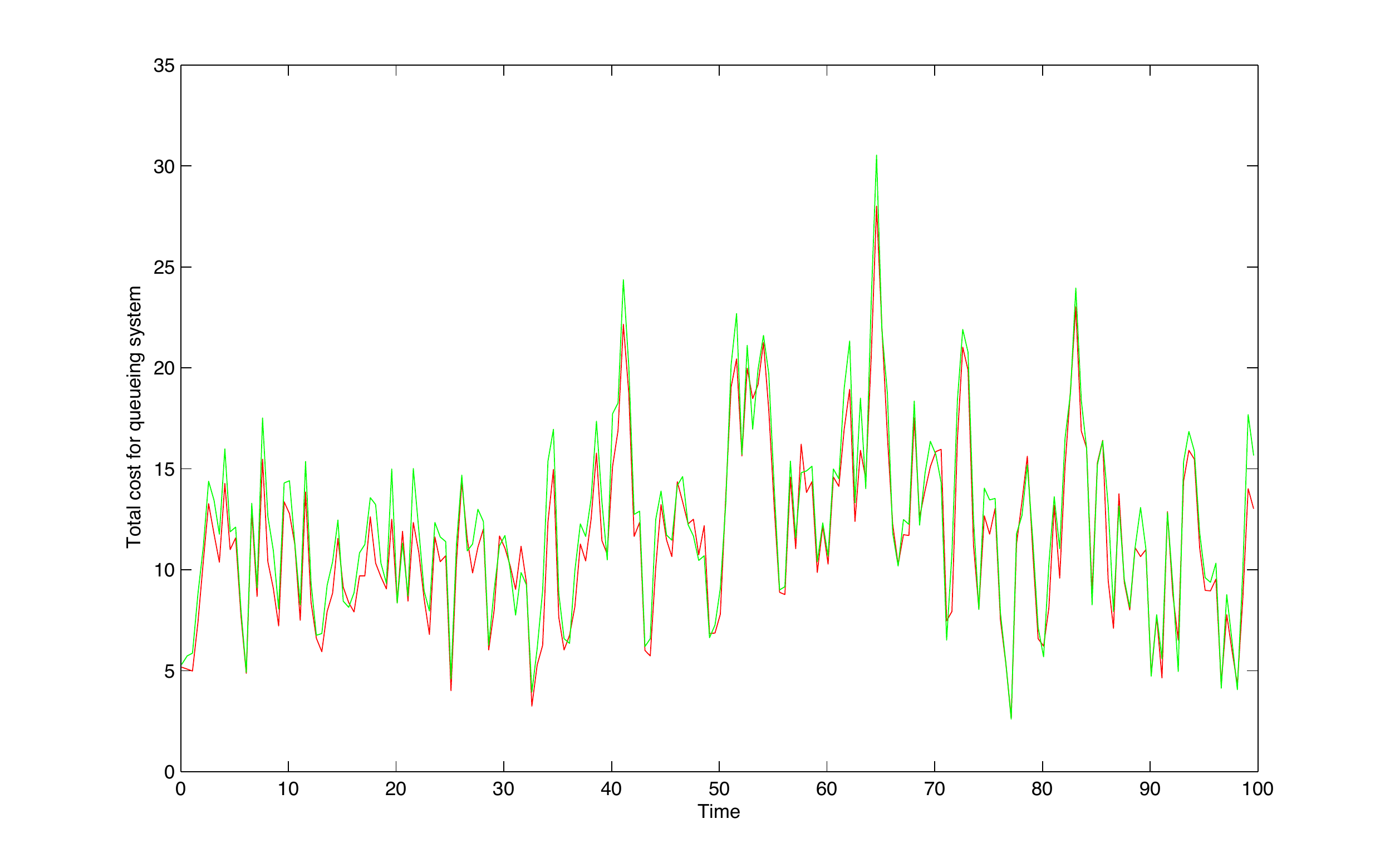}
\includegraphics[width=160mm, height=45mm]{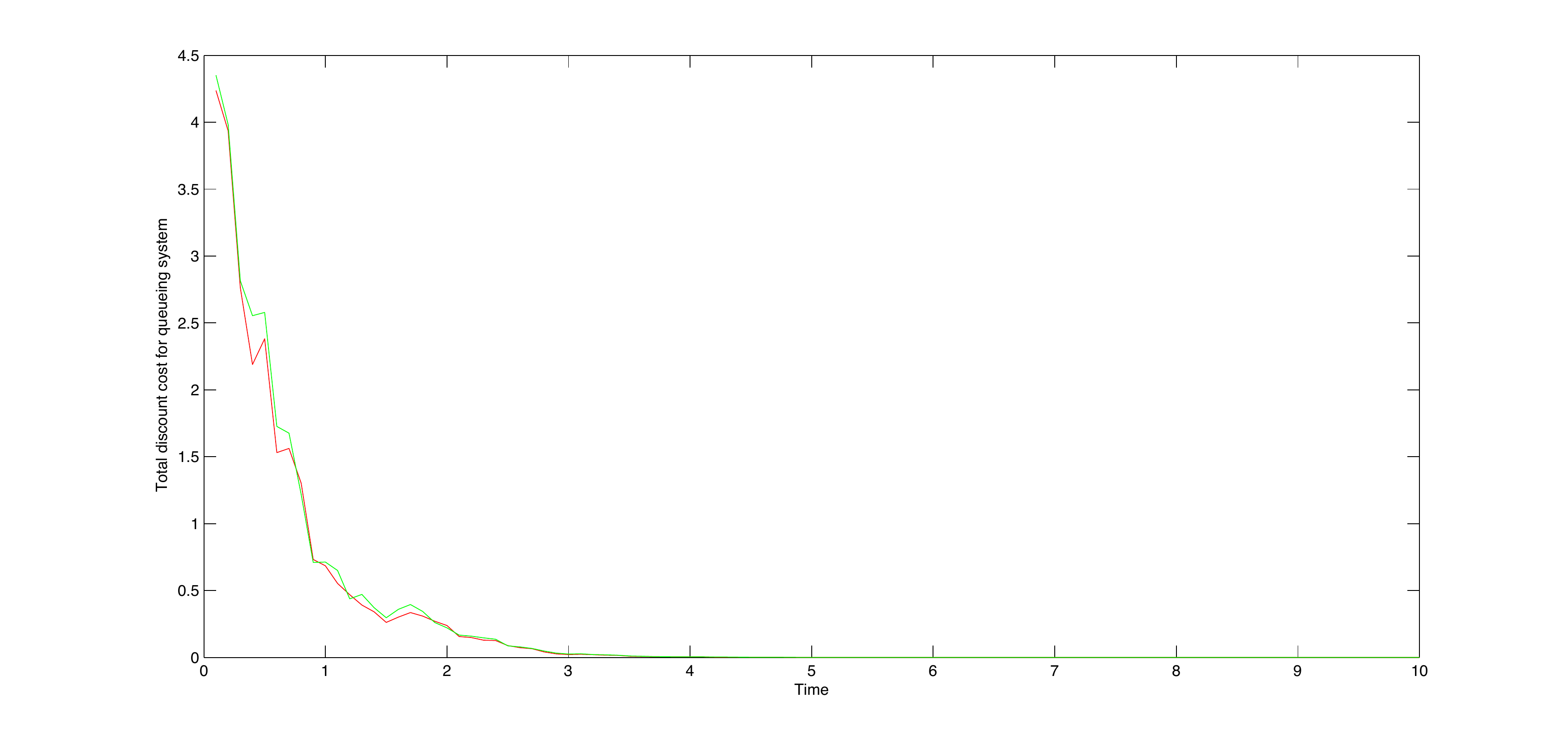}
\caption{Average sample path for the case when $\nu= -1/3$ and $\alpha=2/3$. }
\end{figure}

\bibliographystyle{plain}
\bibliography{reference}

\end{document}